\documentclass[12pt,reqno]{amsart}
\usepackage{amsfonts}
\usepackage{amsthm}
\usepackage{amsmath,amssymb}
\usepackage{amscd}
\usepackage[latin2]{inputenc}
\usepackage{t1enc}
\usepackage[mathscr]{eucal}
\usepackage{indentfirst}
\usepackage{graphicx}
\usepackage{graphics}
\usepackage{pict2e}
\usepackage{epic}
\usepackage[margin=2.9cm]{geometry}
\usepackage{epstopdf} 
\usepackage{bm}
\usepackage[all]{xy}
\usepackage{mathrsfs}
\usepackage{lmodern}
\newtheoremstyle{break}%
{}{}%
{\itshape}{}%
{\bfseries}{}
{\newline}{}

\newtheoremstyle{dbreak}%
{}{}%
{\upshape}{}%
{\bfseries}{}
{\newline}{}

\theoremstyle{break}
\newtheorem{Thm}{Theorem}[section]

\newtheorem{Prop}[Thm]{Proposition}
\newtheorem{Cor}[Thm]{Corollary}

\newtheorem{theorem}{Theorem}

\theoremstyle{dbreak}
\newtheorem{Def}[Thm]{Definition}
\newtheorem{Rem}[Thm]{Remark}

\newcommand{\prf}{\noindent\underline{$Proof.$}}

\numberwithin{equation}{section}

\newenvironment{nouppercase}{%
\renewcommand{\uppercasenonmath}[1]{}}{}

\makeatletter
\@namedef{subjclassname@2020}{%
  \textup{2020} Mathematics Subject Classification}
\makeatother

\DeclareMathOperator*{\Star}{\text{\Huge{$\ast$}}}

\usepackage[hyperindex,breaklinks]{hyperref}
\PassOptionsToPackage{linktocpage}{hyperref}

\begin{document}

\title[Tame fundamental groups of singular curves]{Anabelian properties of tame fundamental groups of singular curves}
\author{Takahiro Murotani}



\address[Takahiro Murotani]{Faculty of Arts and Sciences, Kyoto Institute of Technology, 
Matsugasaki, Sakyo-ku, Kyoto 606-8585, Japan}

\email{murotani@kit.ac.jp}

\begin{abstract}
In anabelian geometry, we consider to what extent the \'{e}tale or tame fundamental groups of schemes reflect geometric properties of the schemes.
Although there are many known results (mainly for smooth curves) in this area, general singular curves have rarely been treated.
One reason is that we cannot determine the isomorphism classes of singular curves themselves from their \'{e}tale or tame fundamental groups in general.
On the other hand, Das proved that the structure of the tame fundamental group of a singular curve over an algebraically closed field is determined by its normalization and an invariant relating to its singularities.
In the present paper, by using this and known anabelian results, we determine the isomorphism classes of the normalizations of singular curves over various fields under certain conditions.
\end{abstract}

\begin{nouppercase}
\maketitle
\end{nouppercase}

\tableofcontents

\section*{Introduction}

In anabelian geometry, we consider to what extent the \'{e}tale or tame fundamental groups (or, more generally, certain quotients of them) of schemes reflect geometric properties of the schemes.
For example, the \'{e}tale fundamental groups of hyperbolic curves over a generalized sub-$p$-adic field $k$ (with the natural surjections to the absolute Galois group of $k$) determine the ($k$-)isomorphism classes of curves (cf. \cite[Theorem 4.2]{sur}). 
(For the definitions of a generalized sub-$p$-adic field and a hyperbolic curve, see \S \ref{NC}.)
Similar results are known for hyperbolic curves over finitely generated fields of positive characteristic (see \cite[Theorem 0.5]{Tama1} and \cite[Theorem 3.12]{cusp} for the case of finite fields, and \cite[Theorem 3.2]{St1} and \cite[Theorem 5.1.1]{St2} for the case of transcendental extensions of prime fields).
These are results for smooth curves, while there are also some results for stable (log) curves (and their admissible fundamental groups) (see, e.g., \cite{closed}, \cite{comb} and \cite{Y}).
However, general singular curves have rarely been treated in this area.
One obvious reason is that we cannot determine the isomorphism classes of singular curves themselves from their \'{e}tale or tame fundamental groups in general (even in the situation where the Grothendieck conjecture for hyperbolic curves is affirmative, cf. Remark \ref{nonisom}).

On the other hand, Das has recently proved that the tame fundamental group of a singular curve over an algebraically closed field is isomorphic to the free product of the tame fundamental groups of the normalizations of the irreducible components of the singular curve and some free profinite group (cf. \cite[Theorem 1.1]{D}).
So it is natural to consider the reconstruction of the isomorphism class of the normalization of a singular curve from its tame fundamental group.
In the present paper, we show that this is possible over various fields (under certain conditions).

We shall introduce some notations to state our main results.
For $\Box\in\{\circ, \bullet\}$, let $k_\Box$ be a field of characteristic $p_{k_\Box}\,(\in\{0\}\cup\mathfrak{Primes})$, $k^\mathrm{s}_\Box$ a separable closure of $k_\Box$, $G_{k_\Box}$ the Galois group $\mathrm{Gal}(k^\mathrm{s}_\Box/k_\Box)$, $X_\Box$ a projective curve over $k_\Box$ (i.e., a geometrically reduced, geometrically connected and separated scheme of finite type over $k_\Box$ of dimension one, but possibly singular or reducible, cf. \S\ref{NC}), $\nu:X^\nu_\Box\to X_\Box$ the normalization map of $X_\Box$, and $\Sigma_\Box$ a nonempty set of prime numbers.
Moreover, we shall denote the maximal pro-$\Sigma_\Box$ quotient of the tame fundamental group $\pi_1^\mathrm{t}((X_\Box)_{k^\mathrm{s}_\Box})$ of $(X_\Box)_{k^\mathrm{s}_\Box}:=X_\Box\otimes_{k_\Box} k^\mathrm{s}_\Box$ (for some choice of basepoint) by $\Delta_\Box$, and set $\Pi_\Box:=\pi_1^\mathrm{t}(X_\Box)/\mathrm{Ker}(\pi_1^\mathrm{t}((X_\Box)_{k^\mathrm{s}_\Box})\twoheadrightarrow\Delta_\Box)$.
Note that, for each $\Box\in\{\circ, \bullet\}$, we have the following exact sequence:

\[\xymatrix{1 \ar[r] & \Delta_\Box \ar[r] & \Pi_\Box \ar[r]^{\mathrm{pr}_\Box} & G_{k_\Box} \ar[r] & 1.}\]
($\Pi_\Box$ surjects onto $G_{k_\Box}$ since $X_\Box$ is geometrically connected.)

For projective singular curves, we have the following:

\begin{theorem}[cf. Remark \ref{sigma} and Theorem \ref{projective}]\label{mainA}
Suppose that, for $\Box\in\{\circ, \bullet\}$, each irreducible component of the normalization of $(X_\Box)_{k^\mathrm{s}_\Box}$ is of genus at least two, and that $\Sigma_\Box$ is an infinite set.
Suppose, moreover, that there exists an isomorphism of profinite groups $\phi:\Pi_\circ\stackrel{\sim}{\to}\Pi_\bullet$ which induces an isomorphism $\Delta_\circ\stackrel{\sim}{\to}\Delta_\bullet$.
Then it holds that $\Sigma_\circ=\Sigma_\bullet$.
(Set 
$\Sigma:=\Sigma_\circ=\Sigma_\bullet$ and $\Sigma':=\mathfrak{Primes}\setminus\Sigma$ (cf. Remark \ref{sigma}).)
Furthermore, the following hold:
\begin{enumerate}
\item[(i)] Suppose that $k_\circ$ and $k_\bullet$ are finite fields.
Suppose, moreover, that $\Sigma'$ is finite.
Then $X_{\circ}^\nu$ and $X_\bullet^\nu$ are isomorphic as schemes.

\item[(ii)] Suppose that $k_\circ=k_\bullet$, and $k:=k_\circ=k_\bullet$ is a finitely generated transcendental extension field of a finite field.
Suppose, moreover, that $p_{k_\circ}=p_{k_\bullet}\not\in\Sigma$, that $\Sigma'$ is finite, that $\phi$ induces the identity automorphism of $G_{k_\circ}=G_{k_\bullet}$, and that any irreducible component of the normalization of $(X_\circ)_{k^\mathrm{s}_\circ}$ is not isotrivial (i.e., not defined over any finite field).
Then, after replacing the irreducible components of $X_\circ^\nu$ and $X_\bullet^\nu$ by suitable Frobenius twists of them, $X_\circ^\nu$ and $X_\bullet^\nu$ are isomorphic as $k$-schemes.

\item[(iii)] Suppose that $k_\circ$ and $k_\bullet$ are the perfections of finitely generated transcendental extensions of finite fields.
Suppose, moreover, that $p_{k_\circ}, p_{k_\bullet}\not\in\Sigma$, that $\Sigma'$ is finite, and that any irreducible component of the normalization of $(X_\circ)_{k^\mathrm{s}_\circ}$ is not isotrivial.
Then $X_\circ^\nu$ and $X_\bullet^\nu$ are isomorphic as schemes.

\item[(iv)] Suppose that one of the following holds:
\begin{enumerate}
\item[(a)] $k_\circ$ and $k_\bullet$ are $p$-adic local fields for some $p\in\mathfrak{Primes}$, and $p\in\Sigma$.

\item[(b)] $k_\circ$ and $k_\bullet$ are finitely generated extensions of $\mathbb{Q}$.
\end{enumerate}
Then $X_{\circ}^\nu$ and $X_\bullet^\nu$ are isomorphic as schemes.

\item[(v)] Suppose that $k_\circ=k_\bullet$, and $k:=k_\circ=k_\bullet$ is a generalized sub-$p$-adic field (for some $p\in\mathfrak{Primes}$).
Suppose, moreover, that $p\in\Sigma$, and that $\phi$ induces the identity automorphism of $G_{k_\circ}=G_{k_\bullet}$.
Then $X_{\circ}^\nu$ and $X_\bullet^\nu$ are isomorphic as $k$-schemes.
\end{enumerate}
\end{theorem}

Moreover, for $\Box\in\{\circ, \bullet\}$, let $D_\Box$ be a (possibly empty) closed subscheme of $X_\Box$ which is finite \'{e}tale over $k_\Box$, and contained in the smooth locus of $X_\Box$, and set $U_\Box:=X_\Box\setminus D_\Box$ and $U^\nu_\Box:=\nu^{-1}(U_\Box)$.
Similarly to the projective case, let $\Delta_{U_\Box}$ be the maximal pro-$\Sigma_\Box$ quotient of the tame fundamental group $\pi_1^\mathrm{t}((U_\Box)_{k^\mathrm{s}_\Box})$ of $(U_\Box)_{k^\mathrm{s}_\Box}:=U_\Box\otimes_{k_\Box}k^\mathrm{s}_\Box$ (for some choice of basepoint), and set $\Pi_{U_\Box}:=\pi_1^\mathrm{t}(U_\Box)/\mathrm{Ker}(\pi_1^\mathrm{t}((U_\Box)_{k^\mathrm{s}_\Box})\twoheadrightarrow\Delta_{U_\Box})$.

The following is a result for singular curves not necessarily projective:

\begin{theorem}[cf. Theorems \ref{mainfin} and \ref{mainsubp}]\label{mainB}
Suppose that, for $\Box\in\{\circ, \bullet\}$, $\Sigma_\Box$ is an infinite set.
Suppose, moreover, that there exists an isomorphism of profinite groups $\phi:\Pi_{U_\circ}\stackrel{\sim}{\to}\Pi_{U_\bullet}$ which induces an isomorphism $\Delta_{U_\circ}\stackrel{\sim}{\to}\Delta_{U_\bullet}$.
Then the following hold:
\begin{enumerate}
\item[(i)] Suppose that, for $\Box\in\{\circ, \bullet\}$, $k_\Box$ is a finite field, and that $\mathfrak{Primes}\setminus\Sigma_\Box$ is a finite set.
Suppose, moreover, that each irreducible component of the inverse image of $(U_\Box)_{k^\mathrm{s}_\Box}$ by the normalization map of $(X_\Box)_{k^\mathrm{s}_\Box}$ is a curve of type distinct from $(0, 0)$, $(0, 1)$ and $(1, 0)$ (i.e., isomorphic to none of $\mathbb{P}^1_{k^\mathrm{s}_\Box}$, $\mathbb{A}^1_{k^\mathrm{s}_\Box}$ and elliptic curves).
Then $U^\nu_\circ$ and $U^\nu_\bullet$ are isomorphic as schemes.

\item[(ii)] Suppose that $k_\circ=k_\bullet$, that $k:=k_\circ=k_\bullet$ is a generalized sub-$p$-adic field (for some $p\in\mathfrak{Primes}$), that $\Sigma_\circ$ and $\Sigma_\bullet$ contain $p$, and that $\phi$ induces the identity automorphism of $G_{k_\circ}=G_{k_\bullet}$.
Suppose, moreover, that each irreducible component of the inverse image of $(U_\Box)_{k^\mathrm{s}_\Box}$ by the normalization map of $(X_\Box)_{k^\mathrm{s}_\Box}$ is a hyperbolic curve.
Then $U^\nu_\circ$ and $U^\nu_\bullet$ are isomorphic as $k$-schemes.
\end{enumerate}
\end{theorem}

The key step of the proofs of these results is the group-theoretic reconstruction of the decomposition groups of the irreducible components of singular curves (from $\Delta_{\Box}\subset\Pi_{\Box}$ or $\Delta_{U_\Box}\subset\Pi_{U_\Box}$).
Roughly speaking, these decomposition groups coincide with the (geometrically pro-$\Sigma_\Box$) tame fundamental groups of the normalizations of the irreducible components (cf. Remarks \ref{decomp} and \ref{afdecomp}).
Then we can apply known anabelian results (for smooth curves) to obtain the main theorems.

The present paper is organized as follows.
In \S \ref{Sproj}, we investigate the structures of the (geometrically pro-$\Sigma_\Box$) tame fundamental groups of projective singular curves over separably closed fields.
In this case, the decomposition groups of the irreducible components (whose normalizations are of positive genus) can be reconstructed group-theoretically without Galois actions.
This is mainly because the number of the irreducible components (whose normalizations are of positive genus) is easily reconstructed group-theoretically.
In \S \ref{various}, we apply results obtained in \S 1 to projective singular curves over various fields, and prove Theorem \ref{mainA}.
In \S \ref{FF} and \S \ref{GENP}, we consider the (geometrically pro-$\Sigma_\Box$) tame fundamental groups of singular curves which are not necessarily projective.
In this case, we can do nothing without Galois actions in general.
(For example, if the base field is of characteristic zero, or if the base field is of characteristic $p(>0)$ and $p\not\in\Sigma_\Box$, then $\Delta_{U_\Box}$ is possibly a free pro-$\Sigma_\Box$ group.)
Moreover, in order to reconstruct the decomposition groups, we need to reconstruct various invariants of singular curves group-theoretically (from $\Delta_{U_\Box}\subset\Pi_{U_\Box}$).
We carry out these reconstructions for singular curves over finite fields (resp. over generalized sub-$p$-adic fields) in \S \ref{FF} (resp. in \S \ref{GENP}), and prove (i) (resp. (ii)) of Theorem \ref{mainB}.

\

\section*{Acknowledgments}

The author would like to thank Professor Akio Tamagawa for various helpful discussions and comments.
The author was supported by JSPS KAKENHI Grant Numbers 22J00022 and 24K16890.

\

\

\setcounter{section}{-1}

\section{Notations and conventions}\label{NC}

\noindent
{\bf{Numbers and fields:}}

\

We shall write
\begin{itemize}
\item $\mathbb{Z}$ for the set of integers;

\item $\mathbb{N}$ for the set of positive integers;

\item $\mathbb{Q}$ for the set of rational numbers;

\item $\mathbb{R}$ for the set of real numbers;

\item $\mathfrak{Primes}$ for the set of prime numbers.

\end{itemize}

For $p\in\mathfrak{Primes}$ and a positive integer $n$, we shall write
\begin{itemize}
\item $\mathbb{Z}_p$ for the $p$-adic completion of $\mathbb{Z}$;

\item $\mathbb{Q}_p$ for the quotient field of $\mathbb{Z}_p$;

\item $\mathbb{F}_{p^n}$ for the finite field of cardinality $p^n$.
\end{itemize}

For a set of prime numbers $\Sigma$, we shall write:
\[N(\Sigma):=\{n\in\mathbb{N}\,|\,\text{if $p\in\mathfrak{Primes}$ divides $n$, then $p\in\Sigma$}\}.\]

We shall call a field $k$ a {\it{generalized sub-$p$-adic field}} if $k$ is isomorphic to a subfield of a finitely generated extension of the quotient field of $W(\overline{\mathbb{F}}_p)$, where $W(\overline{\mathbb{F}}_p)$ is the ring of Witt vectors with coefficients in an algebraic closure $\overline{\mathbb{F}}_p$ of $\mathbb{F}_p$.

\

\noindent
{\bf{Sets:}}

\

For a finite set $X$, we shall write $|X|$ for the cardinality of $X$.

\

\noindent
{\bf{Profinite groups:}}

\

Let $G$ be a profinite group, $p\in\mathfrak{Primes}$ and $\Sigma$ a nonempty set of prime numbers.
Then we shall write $G^\mathrm{ab}$ for the abelianization of $G$ (i.e., the quotient of $G$ by the closure of the commutator subgroup of $G$), $G^p$ (resp. $G^{p'}$, resp. $G^{\Sigma}$) for the maximal pro-$p$ quotient (resp. the maximal pro-prime-to-$p$ quotient, resp. the maximal pro-$\Sigma$ quotient) of $G$, and $G^{p\text{-ab}}$ for the abelianization of $G^p$.
We denote the cohomological $p$-dimension of $G$ by $\mathrm{cd}_pG$.
Moreover, for two profinite groups $H_1$ and $H_2,$ we denote the free profinite product of $H_1$ and $H_2$ by $H_1\ast H_2$.

\

\noindent
{\bf{Curves:}}

\

Let $U$ be a scheme over a field $k$.
We shall say that $U$ is a {\it{curve over $k$}} if $U$ is a geometrically reduced, geometrically connected and separated scheme of finite type over $k$ of dimension one.

Moreover, suppose that there exist a projective curve $X$ over $k$ and a (possibly empty) closed subscheme $D$ of $X$ which is finite \'{e}tale over $k$ and contained in the smooth locus of $X$ satisfying $U=X\setminus D$.
Then we denote the tame fundamental group of $X$ with respect to $D$ (for some choice of base point) by $\pi_1^\mathrm{t}(U)$.
Suppose, furthermore, that $X$ is smooth (hence geometrically irreducible).
Let $g$ be the genus of $X$ and $r:=|D(k^\mathrm{s})|$, where $k^\mathrm{s}$ is a separable closure of $k$.
Then we shall say that $U$ is a {\it{curve of type}} $(g, r)$ ({\it{over}} $k$).
Moreover, if $2g+r-2>0$, we shall say that $U$ is a {\it{hyperbolic curve}} ({\it{over}} $k$).

\

\

\section{Fundamental groups and irreducible components of projective singular curves over separably closed fields}\label{Sproj}

Let $k$ be a separably closed field of characteristic $p_k\,(\in\{0\}\cup\mathfrak{Primes})$.
Moreover, we shall write:

\begin{itemize}
\item $X$ for a projective curve over $k$ (possibly singular or reducible);

\item $n$ for the number of irreducible components of $X$;

\item $\nu: X^\nu\to X$ for the normalization map;

\item $X_1, \cdots, X_n\subset X$ for the (distinct) irreducible components of $X$;

\item $X_i^\nu$ for the normalization of $X_i\, (1\leq i\leq n)$;

\item $g_i$ for the genus of $X_i^\nu\, (1\leq i\leq n)$;

\item $n^{>0}:=|\{i\in\mathbb{Z}\,|\,1\leq i\leq n, g_i>0\}|$;

\item $\Sigma\subset\mathfrak{Primes}$ for a nonempty set of prime numbers;

\item $\Sigma^\dagger:=\Sigma\setminus\{p_k\}$;

\item $\Delta_X:=\pi^\mathrm{t}_1(X)^\Sigma$, $\Delta_{X_i}:=\pi^\mathrm{t}_1(X_i)^\Sigma$, $\Delta_{X_i^\nu}:=\pi_1^\mathrm{t}(X_i^\nu)^\Sigma$.
\end{itemize}

Furthermore, if $p_k>0$, for $1\leq i\leq n$, let $\gamma_i$ be the $p_k$-rank of the Jacobian variety of $X_i^\nu$.

By \cite[Theorem 1.1]{D} (see also \cite[Theorem 1.3]{Hove}), we have the following isomorphism of profinite groups:
\begin{align}\label{isom}
\Delta_X\simeq\left(\Star_{i=1}^n\Delta_{X_i^\nu}\right)\ast F_\delta,
\end{align}
where $\delta:=1-n+\displaystyle\sum_{x\in X}\left(\left|\nu^{-1}(x)\right|-1\right)$ and $F_\delta$ is a free pro-$\Sigma$ group of rank $\delta$.

In the following, we suppose that $|\Sigma|\geq 2$.

\begin{Rem}\label{p}
If $n^{>0}\neq 0$ (resp. $\delta\neq 0$), we may determine $\Sigma^\dagger$ (resp. $\Sigma$) group-theoretically from $\Delta_X$.
Indeed, if $n^{>0}\neq 0$, the following conditions are equivalent for $l\in\mathfrak{Primes}$:
\begin{enumerate}
\item[(i)] $l\in\Sigma^\dagger$.

\item[(ii)] $\Delta_X^l$ is not trivial, and for any $l'\in\mathfrak{Primes}$, it holds that $\mathrm{rank}_{\mathbb{Z}_l}\Delta_X^{l\text{-ab}}\geq\mathrm{rank}_{\mathbb{Z}_{l'}}\Delta_X^{l'\text{-ab}}$.
\end{enumerate}
On the other hand, if $\delta\neq 0$, then we have $\Sigma=\{l\in\mathfrak{Primes}\,|\,\Delta_X^l\text{ is not trivial}\}$.

Moreover, if one of the following conditions holds, we may also determine group-theoretically whether $\Sigma$ contains $p_k$ or not:
\begin{enumerate}
\item[(a)] $g_i\geq 2$ for some $1\leq i\leq n$.

\item[(b)] $n^{>0}\neq 0$ and $\delta\neq 0$.
\end{enumerate}
Indeed, in this case, by \cite[Lemma 1.9]{Tama1}, $p_k\in\Sigma$ if and only if there exist an open subgroup $\mathcal{H}$ of $\Delta_X$ and $l\in\mathfrak{Primes}$ such that $\mathcal{H}^l$ is a nontrivial free pro-$l$ group (note that $|\Sigma|\geq 2$).
Furthermore, if the latter condition holds, then it holds that $l=p_k$.

\end{Rem}

By renumbering if necessary, we assume that $g_i\geq 1$ (resp. $g_i=0$) for $1\leq i\leq n^{>0}$ (resp. $n^{>0}+1\leq i\leq n$).
Set $\Delta_i:=\Delta_{X_i^\nu}$ for $1\leq i\leq n$ and $\displaystyle \Delta:=\left(\Star_{i=1}^n\Delta_i\right)\ast F_\delta=\left(\Star_{i=1}^{n^{>0}}\Delta_i\right)\ast F_\delta$, and fix an isomorphism $\iota:\Delta_X\stackrel{\sim}{\to}\Delta$ (cf. (\ref{isom})).
Moreover, for an open subgroup $\mathcal{H}\subset\Delta_X$, we denote the covering of $X$ corresponding to $\mathcal{H}$ by $X_\mathcal{H}$, and invariants of $X_\mathcal{H}$ by the corresponding invariants of $X$ with subscript $\mathcal{H}$ (e.g., $n_\mathcal{H}$, $(n^{>0})_\mathcal{H}$ etc.)

\begin{Prop}\label{delta}
$n^{>0}$ and $\delta$ are group-theoretically recovered from $\Delta_X$.
\end{Prop}

\prf

If $\Delta_X$ is trivial, then it holds that $n^{>0}=\delta=0$.
So we assume that $\Delta_X$ is not trivial.
First, we shall reconstruct $n^{>0}$.
For $l\in\mathfrak{Primes}$, by \cite[Proposition 4.1.3 and Theorem 4.1.4]{NSW},
\[H^2(\Delta_X^l, \mathbb{Z}/l\mathbb{Z})\simeq \left(\bigoplus_{i=1}^nH^2(\Delta_i^l, \mathbb{Z}/l\mathbb{Z})\right)\oplus H^2(F_\delta^l, \mathbb{Z}/l\mathbb{Z}).\]
Therefore, $n^{>0}=\displaystyle\max_{l\in\mathfrak{Primes}}\dim_{\mathbb{Z}/l\mathbb{Z}}H^2(\Delta_X^l, \mathbb{Z}/l\mathbb{Z})$ is group-theoretically recovered from $\Delta_X$ (note that $|\Sigma|\geq 2$).
If $n^{>0}=0$, then we have $\delta=\displaystyle\max_{l\in\mathfrak{Primes}}\mathrm{rank}_{\mathbb{Z}_l}\Delta_X^{l\text{-ab}}$.
In the following, we assume that $n^{>0}\neq 0$

We shall reconstruct $\delta$.
By replacing $\Delta_X$ by its maximal pro-$\Sigma^\dagger$ quotient if necessary, we assume that $p_k\not\in\Sigma$ (and $|\Sigma|\geq 1$) (cf. Remark \ref{p}).
Let $\mathcal{H}\subset\Delta_X$ be an open normal subgroup of index $N\,(\in\mathbb{N}(\Sigma))$, and set $H:=\iota(\mathcal{H})$.

First, consider the case where $\delta>0$.
Then, by \cite[Theorem 9.1.9]{RZ}, we have:
\[H=\left(\Star_{i=1}^{n^{>0}}\left(\Star_{\tau\in H\backslash \Delta/\Delta_i}(H\cap h_{i, \tau}\Delta_ih_{i, \tau}^{-1})\right)\right)\ast\left(\Star_{\tau\in H\backslash \Delta/F_{\delta}}(H\cap h_{\tau}F_\delta h_{\tau}^{-1})\right)\ast F,\]
where $h_{i, \tau}$ (resp. $h_\tau$) ranges over a system of double coset representations for $H\backslash \Delta/\Delta_i$ (resp. $H\backslash \Delta/F_\delta$) containing the unit element, and $F$ is a free pro-$\Sigma$ group of rank
\[\displaystyle 1+n^{>0}N-\sum_{i=1}^{n^{>0}}|H\backslash \Delta/\Delta_i|-|H\backslash \Delta/F_\delta|.\]
For each $1\leq i\leq n^{>0}$ and $\tau\in H\backslash \Delta/\Delta_i$, $|H\backslash \Delta/\Delta_i|=|\Delta/H\Delta_i|$ and $|h_{i, \tau}\Delta_ih_{i, \tau}^{-1}/H\cap h_{i, \tau}\Delta_ih_{i, \tau}^{-1}|=|\Delta_i/\Delta_i\cap H|=|H\Delta_i/H|$ (note that $H$ is an open normal subgroup of $\Delta$), and similar formulas hold also for $F_\delta$.
Therefore, by the Hurwitz formula, for each $1\leq i\leq n^{>0}$ and $\tau\in H\backslash \Delta/\Delta_i$, we have
\[\mathrm{rank}_\Sigma (H\cap h_{i, \tau}\Delta_i h_{i, \tau}^{-1})^{\mathrm{ab}}=2(|H\Delta_i/H|(g_i-1)+1),\]
where we write the rank of finitely generated free pro-$\hat{\mathbb{Z}}^\Sigma$ module $A$ for $\mathrm{rank}_\Sigma A$.
Similarly, for each $\tau\in H\backslash \Delta/F_\delta$, it holds that
\[\mathrm{rank}_\Sigma(H\cap h_\tau F_\delta h_\tau^{-1})^{\mathrm{ab}}=|HF_\delta/H|(\delta-1)+1.\]
Therefore,
\begin{align*}
&\mathrm{rank}_\Sigma H^{\mathrm{ab}} \\
=&\sum_{i=1}^{n^{>0}}|\Delta/H\Delta_i|\cdot 2\left(|H\Delta_i/H|(g_i-1)+1\right)+|\Delta/HF_\delta|(|HF_\delta/H|(\delta-1)+1) \\
&\hspace{250pt} +1+n^{>0}N-\sum_{i=1}^{n^{>0}}|\Delta/H\Delta_i|-|\Delta/HF_\delta| \\
=&\sum_{i=1}^{n^{>0}} 2|\Delta/H\Delta_i|\cdot |H\Delta_i/H|(g_i-1)+\sum_{i=1}^{n^{>0}} |\Delta/H\Delta_i|+|\Delta/HF_\delta|\cdot |HF_\delta/H|(\delta-1)+1+n^{>0}N \\
=&2N\sum_{i=1}^{n^{>0}}(g_i-1)+\sum_{i=1}^{n^{>0}}|\Delta/H\Delta_i|+N(\delta-1+n^{>0})+1.
\end{align*}
A similar calculation shows that, for any open normal subgroup $\mathcal{H}\subset\Delta_X$ of index $N$ (and $H:=\iota(\mathcal{H})$), the following holds even if $\delta=0$:
\[\mathrm{rank}_\Sigma H^{\mathrm{ab}}=2N\sum_{i=1}^{n^{>0}}(g_i-1)+\sum_{i=1}^{n^{>0}} |\Delta/H\Delta_i|+N(\delta-1+n^{>0})+1.\]

For $N\in\mathbb{N}(\Sigma)$, set:
\begin{align*}
M(N)&:=\max\{\mathrm{rank}_\Sigma H^{\mathrm{ab}}\,|\,H\subset \Delta:\text{ an open normal subgroup}, (\Delta:H)=N\}, \\
m(N)&:=\min\{\mathrm{rank}_\Sigma H^{\mathrm{ab}}\,|\,H\subset \Delta:\text{ an open normal subgroup}, (\Delta:H)=N\}.
\end{align*}
Then, for $l\in\Sigma$, we have
\begin{align}\label{M}
M(l)-m(l)=\begin{cases} (n^{>0}-1)(l-1)\, &(\delta=0), \\ n^{>0}(l-1)\, &(\delta>0).\end{cases}
\end{align}
Note that, for an open normal subgroup $H\subset \Delta$ satisfying $(\Delta:H)=l$, $H$ satisfies $\mathrm{rank}_\Sigma H^{\mathrm{ab}}=m(l)$  if and only if, $f_i$ is nontrivial for $1\leq i\leq n^{>0}$, where $f=(f_i)_{i=1}^{n^{>0}}+f_\delta\in\displaystyle\left(\bigoplus_{i=1}^{n^{>0}} H^1(\Delta_i, \mathbb{Z}/l\mathbb{Z})\right)\oplus H^1(F_\delta, \mathbb{Z}/l\mathbb{Z})=H^1(\Delta, \mathbb{Z}/l\mathbb{Z})$ is the element corresponding to the natural surjection $\Delta\twoheadrightarrow \Delta/H$.
Moreover, $H$ satisfies $\mathrm{rank}_\Sigma H^{\mathrm{ab}}=M(l)$ if and only if one (and only one) of the following conditions holds:
\begin{enumerate}
\item[(i)] $\delta=0$, and there exists $1\leq j\leq n^{>0}$ such that, $f_j$ is nontrivial, and that $f_i$ is trivial for $1\leq i\leq n^{>0}, i\neq j$.

\item[(ii)] $\delta>0$, and $f_i$ is trivial for $1\leq i\leq n^{>0}$.
\end{enumerate}

Since $n^{>0}$ is group-theoretic as above, we may determine group-theoretically whether $\delta=0$ or not. 
In the following, we assume that $\delta>0$.
Set:
\[\mathfrak{\Delta}:=\{H\subset \Delta:\text{ an open normal subgroup}\,|\,\mathrm{rank}_\Sigma H^{\mathrm{ab}}=M((\Delta:H))\}.\]
Then, clearly, the subgroup $\displaystyle H_\mathfrak{\Delta}:=\bigcap_{H\in\mathfrak{\Delta}}H$ (of $\Delta$) coincides with the smallest closed normal subgroup of $\Delta$ containing $\Delta_1, \cdots, \Delta_{n^{>0}}$.
Therefore, by \cite[IV, \S 1]{NSW},
\[\Delta/H_\mathfrak{\Delta}\simeq F_\delta,\]
and hence $\delta=\mathrm{rank}_\Sigma(\Delta/H_\mathfrak{\Delta})^{\mathrm{ab}}$.
\qed

\

Until the end of this section, we suppose that $n^{>0}\neq 0$.
By a similar argument to \cite[\S 1]{closed}, we may recover the set of irreducible components of $X$ whose normalizations are of genus at least one group-theoretically from $\Delta_X$:

\begin{Prop}\label{irred}
For any open normal subgroup $\mathcal{H}\subset\Delta_X$, the set $I_{X_{\mathcal{H}}}$ of irreducible components of $X_{\mathcal{H}}$ whose normalizations are of genus at least one, together with the natural action of $\Delta_X/\mathcal{H}$ on $I_{X_\mathcal{H}}$, is group-theoretically recovered from $\mathcal{H}\subset\Delta_X$.
Moreover, if $\Sigma$ is an infinite set, for any open normal subgroup $\mathcal{H}\subset\Delta_X$, the natural map $I_{X_{\mathcal{H}}}\to I_X$ is also group-theoretically recovered from $\mathcal{H}\subset\Delta_X$, where $I_X=I_{X_{\Delta_X}}$.
\end{Prop}

\prf

First, we prove the first assertion for $\mathcal{H}=\Delta_X$.
Fix a prime number $l\in\Sigma^\dagger$.
Let $\pi:\Delta\twoheadrightarrow F_\delta$ be the natural surjection obtained in the proof of Proposition \ref{delta} (in a group-theoretic way), and set $L^c:=H^1(F_\delta, \mathbb{Z}/l\mathbb{Z}), L^e:=H^1(\Delta, \mathbb{Z}/l\mathbb{Z})$.
Then $\pi$ induces an injection $L^c\hookrightarrow L^e$, and let $L^n$ be the cokernel of this injection, i.e., the following sequence is exact:
\[\xymatrix{0 \ar[r] & L^c \ar[r] & L^e \ar[r] & L^n \ar[r] & 0.}\]
Note that, as in \cite[\S 1]{closed}, we have $L^n=\displaystyle\bigoplus_{i=1}^n H^1(\Delta_i, \mathbb{Z}/l\mathbb{Z})\simeq \bigoplus_{Z\in I_X}H^1_{\'{e}t}(Z, \mathbb{Z}/l\mathbb{Z})$.

Elements of $L^e$ correspond to tame (hence also \'{e}tale) coverings of $X$ of degree $l$.
Let $L^\ast\subset L^e$ be the subset of elements of $L^e$ whose images in $L^n$ are nonzero.
Let us take an element $\alpha\in L^\ast$.
Let $Y_\alpha\to X$ be the corresponding covering, $\nu_\alpha: Y_\alpha^\nu\to Y_\alpha$ the normalization map, $n^{>0}_\alpha$ the number of irreducible components of $Y_\alpha$ whose normalizations are of genus at least one.
Then, by sending $\alpha\in L^\ast$ to $n^{>0}_\alpha\in\mathbb{Z}$, we obtain a map $\varepsilon:L^\ast\to\mathbb{Z}$.
Since $L^\ast$ is a finite set, we have $m:=\displaystyle\max_{\alpha\in L^\ast}\varepsilon(\alpha)<\infty$, and set $M:=\varepsilon^{-1}(m)$.
Note that these reconstructions of $L^\ast$, $\varepsilon$ and $M\subset L^\ast$ are all group-theoretic (cf. Proposition \ref{delta}).
As in \cite[\S 1]{closed}, let us define a pre-equivalence relation ``$\sim$'' on $M$ as follows:

\begin{quote}
For $\alpha, \beta\in M$, we write $\alpha\sim\beta$ if, for every $\lambda, \mu\in(\mathbb{Z}/l\mathbb{Z})^\times$ such that $\lambda\cdot\alpha+\mu\cdot\beta\in L^\ast$, we have $\lambda\cdot\alpha+\mu\cdot\beta\in M$.
\end{quote}
Then, as in \cite[Proposition 1.3]{closed}, ``$\sim$'' is an equivalence relation, and $C_X:=M/\sim$ is naturally isomorphic to the set of irreducible components of $X$ whose normalizations are of genus at least one.
This completes the proof of the first assertion of this proposition for $\mathcal{H}=\Delta_X$.
In the general case, the action of $\Delta_X/\mathcal{H}$ on $I_{X_{\mathcal{H}}}$ is recovered by considering the action of $\Delta_X/\mathcal{H}$ on $H^1(\iota(\mathcal{H}), \mathbb{Z}/l\mathbb{Z})$ (via $\iota$).

In the remainder of this proof, suppose that $\Sigma$ is an infinite set.
Let us take an open subgroup $\mathcal{H}\subset\Delta_X$.
Moreover, suppose that $l\in\Sigma^\dagger$ does not divide $(\Delta_X:\mathcal{H})$ (note that $\Sigma$ is infinite).
Let $X':=X_{\mathcal{H}}$ be the finite \'{e}tale covering of $X$ corresponding to $\mathcal{H}$, $\nu':(X')^\nu\to X'$ be the normalization map,
$n_i'$ the number of irreducible components of $X'$ above $X_i$ ($1\leq i\leq n^{>0}$), $X'_{i1}, \cdots, X'_{in_i'}$ the (distinct) irreducible components of $X'$ above $X_i$ ($1\leq i\leq n^{>0}$), and $\delta':=\delta_{\mathcal{H}}=1-n_{\mathcal{H}}+\displaystyle\sum_{x'\in X'}\left(\left|(\nu')^{-1}(x')\right|-1\right)$.
Note that, for any $1\leq i\leq n^{>0}$ and $1\leq j\leq n_i'$, the normalization of $X_{ij}'$ is of genus at least $1$.
Similarly to the case of $X$, we define $H_{ij}:=\Delta_{(X_{ij}')^{\nu'}}$ ($1\leq i\leq n^{>0}, 1\leq j\leq n_i'$).
Then, by (\ref{isom}) and \cite[Theorem 9.1.9]{RZ}, we have the following isomorphism:
\[\mathcal{H}\simeq\iota(\mathcal{H})=\left(\Star_{i=1}^{n^{>0}}\Star_{j=1}^{n_i'}H_{ij}\right)\ast F_{\delta'}=:H,\]
where $F_{\delta'}$ is a free pro-$\Sigma$ group of rank $\delta'$.
By the proof of Proposition \ref{delta}, we can recover the natural surjections $\pi: \Delta\twoheadrightarrow F_\delta$ and $\pi': H\twoheadrightarrow F_{\delta'}$.
These induce homomorphisms $\pi^{\mathrm{ab}}: \Delta^{\mathrm{ab}}\twoheadrightarrow F_\delta^{\mathrm{ab}}$ and $(\pi')^{\mathrm{ab}}: H^{\mathrm{ab}}\twoheadrightarrow F_{\delta'}^{\mathrm{ab}}$.
Moreover, the kernels of these homomorphisms are isomorphic to $\displaystyle\bigoplus_{i=1}^{n^{>0}} \Delta_i^{\mathrm{ab}}$ and $\displaystyle\bigoplus_{i=1}^{n^{>0}}\bigoplus_{j=1}^{n_i'}H_{ij}^{\mathrm{ab}}$, respectively.
We identify the kernels of $\pi^{\mathrm{ab}}$ and $(\pi')^{\mathrm{ab}}$ with these groups.
Therefore, we obtain the following diagram:
\[\xymatrix@C=40pt{\displaystyle 0 \ar[r] & \displaystyle\bigoplus_{i=1}^{n^{>0}}\bigoplus_{j=1}^{n_i'}H_{ij}^{\mathrm{ab}} \ar[r] & \displaystyle H^{\mathrm{ab}}\simeq\left(\bigoplus_{i=1}^{n^{>0}}\bigoplus_{j=1}^{n_i'} H_{ij}^{\mathrm{ab}}\right)\oplus F_{\delta'}^{\mathrm{ab}} \ar[r]^{\hspace{60pt}(\pi')^{\mathrm{ab}}} \ar[d] & F_{\delta'}^{\mathrm{ab}} \ar[r] & 0 \\ 0 \ar[r] & \displaystyle\bigoplus_{i=1}^{n^{>0}}\Delta_i^{\mathrm{ab}} \ar[r] & \displaystyle \Delta^{\mathrm{ab}}\simeq\left(\bigoplus_{i=1}^{n^{>0}} \Delta_i^{\mathrm{ab}}\right)\oplus F_\delta^{\mathrm{ab}} \ar[r]^{\hspace{50pt}\pi^{\mathrm{ab}}} & F_\delta^{\mathrm{ab}} \ar[r] & 0,}\]
where the horizontal sequences are split exact, and the vertical arrow is the natural homomorphism.
Since, for each $1\leq i\leq n^{>0}$ and $1\leq j\leq n_i'$, the image of $H_{ij}^{\mathrm{ab}}$ is contained in $\Delta_i^{\mathrm{ab}}$, the image of any element of $H^1(F_\delta^{\mathrm{ab}}, \mathbb{Z}/l\mathbb{Z})$ in $H^1(\mathrm{Ker}(\pi')^{\mathrm{ab}}, \mathbb{Z}/l\mathbb{Z})$ vanishes.
Therefore, by considering $H^1(-, \mathbb{Z}/l\mathbb{Z})$ of the above diagram, we obtain the following commutative diagram:
\[\xymatrix{ 0 \ar[r] & L^c_H \ar[r] & L^e_H \ar[r] & L^n_H \ar[r] & 0 \\ 0 \ar[r] & L^c \ar@{.>}[u] \ar[r] & L^e \ar[u] \ar[r] & L^n \ar@{.>}[u]_{\phi_H} \ar[r] & 0,}\]
where $L_H^c:=H^1(F_{\delta'}, \mathbb{Z}/l\mathbb{Z})=H^1(F_{\delta'}^{\mathrm{ab}}, \mathbb{Z}/l\mathbb{Z}), L_H^e:=H^1(H, \mathbb{Z}/l\mathbb{Z})=H^1(H^{\mathrm{ab}}, \mathbb{Z}/l\mathbb{Z})$, $L_H^n$ is the cokernel of the injection $L_H^c\hookrightarrow L_H^e$, the horizontal sequences are split exact and the dotted vertical arrows are defined by the above arguments.
Note that we already know $I_X$ and $I_{X'}\,(:=I_{X_\mathcal{H}})$ (in a group-theoretic way).
As in \cite[\S 1]{closed}, we have
\begin{align*}
L^n&=\bigoplus_{i=1}^{n^{>0}} H^1(\Delta_i, \mathbb{Z}/l\mathbb{Z})\simeq\bigoplus_{Z\in I_X}L_Z, \\
L^n_H&=\bigoplus_{i=1}^{n^{>0}}\bigoplus_{j=1}^{n_i'} H^1(H_{ij}, \mathbb{Z}/l\mathbb{Z})\simeq\bigoplus_{Z'\in I_{X'}}L_{Z'},
\end{align*}
where $L_Z:=H^1_{\'{e}t}(Z, \mathbb{Z}/l\mathbb{Z})$ and $L_{Z'}:=H^1_{\'{e}t}(Z', \mathbb{Z}/l\mathbb{Z})$ for each $Z\in I_X$ and $Z'\in I_{X'}$.
Therefore, the following two conditions are equivalent (note that $l$ is prime to $(\Delta_X:\mathcal{H})$):
\begin{enumerate}
\item[(i)] $Z'\in I_{X'}$ lies above $Z\in I_X$.

\item[(ii)] For any nonzero element $\alpha\in L_Z\subset L^n$, the $Z'$-component of $\phi_H(\alpha)\in L_H^n$ is not zero.
\end{enumerate}
On the other hand, we can recover group-theoretically not only $L^n$ and $L_H^n$, but also subspaces corresponding to $L_Z$ and $L_{Z'}$ (for each $Z\in I_X$ and $Z'\in I_{X'}$) by considering the image of each equivalence class of ``$M$'' in $L^n$ (or $L_H^n$) (cf. the proof of the first assertion).
So, the above condition (ii) is a group-theoretic condition.
Note that, as in \cite[\S 1, Remark]{closed}, these constructions do not depend on the choice of $l$.
This completes the proof of the second assertion of this proposition, hence also that of this proposition.
\qed

\begin{Cor}\label{comp}
Suppose that $\Sigma$ is an infinite set.
Then, for each $i$ ($1\leq i\leq n^{>0}$), the natural (up to inner automorphisms) inclusion $\Delta_{X_i^\nu}\hookrightarrow\Delta_X$ is group-theoretically recovered from $\Delta_X$.
In particular, $g_i$ $(1\leq i\leq n^{>0})$ is group-theoretically recovered from $\Delta_X$.
Moreover, if $p_k\in\Sigma$, $\gamma_i$ $(1\leq i\leq n^{>0})$ is also group-theoretically recovered.
\end{Cor}

\prf

Let us fix an integer $i$ $(1\leq i\leq n^{>0})$.
For each open normal subgroup $\mathcal{H}\subset\Delta_X$, let $X_{\mathcal{H}}$ be the finite covering of $X$ corresponding to $\mathcal{H}$, $I_{X_{\mathcal{H}}}$ the set of irreducible components of $X_{\mathcal{H}}$ whose normalizations are of genus at least one and $\psi_{\mathcal{H}}:I_{X_{\mathcal{H}}}\to I_X\,(:=I_{X_{\Delta_X}})$ be the natural map (obtained group-theoretically in Proposition \ref{irred}).
These data determine a projective limit:
\[\mathfrak{I}:=\varprojlim_{\mathcal{H}\subset\Delta_X}I_{X_{\mathcal{H}}},\]
where $\mathcal{H}$ runs through the set of open normal subgroups of $\Delta_X$ and the transition maps are the natural maps $\psi_{\mathcal{H}}^{\mathcal{H'}}:I_{X_{\mathcal{H}'}}\to I_{X_\mathcal{H}}$ for open normal subgroups $\mathcal{H}'\subset\mathcal{H}\subset\Delta_X$.
Note that, by Proposition \ref{irred}, $\mathfrak{I}$, together with the natural action of $\Delta_X$ on $\mathfrak{I}$ can be recovered group-theoretically from $\Delta_X$.
Let $J\in\mathfrak{I}$ be an element which maps to $X_i\in I_X$ and $D_i$ the stabilizer of $J$ (with respect to the action of $\Delta_X$).
Then the conjugacy class of the subgroup $D_i$ of $\Delta_X$ does not depend on the choice of $J$, and the natural inclusion $D_i\hookrightarrow\Delta_X$ coincides with $\Delta_{X_i^\nu}\hookrightarrow\Delta_X$ up to inner automorphisms.
\qed

\

By combining this result with \cite[Theorem 0.1]{Tama4}, we obtain the following result:

\begin{Thm}
Let $k$ be an algebraic closure of $\mathbb{F}_p$ ($p\in\mathfrak{Primes}$), $X$ a (reduced, connected and) projective curve over $k$ (possibly singular or reducible), and $\Delta_X$ the (tame) fundamental group of $X$ (for some choice of basepoint) (i.e., $\Sigma=\mathfrak{Primes}$).
Suppose that the genus of the normalization of each irreducible component of $X$ is at least two.
Then (the isomorphism class of) $\Delta_X$ determines the $k$-isomorphism class of the normalization of $X$ up to finite possibilities.
\end{Thm}

\

\

\section{Anabelian results for projective singular curves over various fields}\label{various}

For $\Box\in\{\circ, \bullet, \emptyset\}$, let $k_\Box$ be a field of characteristic $p_{k_\Box}\,(\in\{0\}\cup\mathfrak{Primes})$.
Moreover, we shall write:

\begin{itemize}

\item $k^\mathrm{s}_\Box$ for a separable closure of $k_\Box$;

\item $G_{k_\Box}$ for the Galois group $\mathrm{Gal}(k^\mathrm{s}_\Box/k_\Box)$;

\item $X_\Box$ for a projective curve over $k_\Box$ (possibly singular or reducible);

\item $n_\Box$ for the number of (distinct) irreducible components of $\overline{X_\Box}:=X_\Box\otimes_{k_\Box}k^\mathrm{s}_\Box$;

\item $\nu_\Box: X^\nu_\Box\to X_\Box$ for the normalization map;

\item $\Sigma_\Box\subset\mathfrak{Primes}$ for a nonempty set of prime numbers;

\item $\Sigma^\dagger_\Box:=\Sigma_\Box\setminus\{p_{k_\Box}\}$;

\item $\Sigma'_\Box:=\mathfrak{Primes}\setminus\Sigma_\Box$;

\item $\Delta_\Box:=\pi_1^\mathrm{t}(\overline{X_\Box})^{\Sigma_\Box}$;

\item $\Pi_\Box:=\pi_1^\mathrm{t}(X_\Box)/\mathrm{Ker}(\pi_1^\mathrm{t}(\overline{X_\Box})\twoheadrightarrow\Delta_\Box)$.
\end{itemize}

Note that, for each $\Box\in\{\circ, \bullet, \emptyset\}$, we have the following exact sequence:

\[\xymatrix{1 \ar[r] & \Delta_\Box \ar[r] & \Pi_\Box \ar[r]^{\mathrm{pr}_\Box} & G_{k_\Box} \ar[r] & 1.}\]
($\Pi_\Box$ surjects onto $G_{k_\Box}$ since $X_\Box$ is geometrically connected.)
Suppose that we are given not only $\Pi_\Box$ but also the subgroup $\Delta_\Box$ as the input datum of various group-theoretic reconstructions. In other words, we are given the above exact sequence as the input data.

In the following, we suppose that each irreducible component of the normalization of $X\otimes_k k^\mathrm{s}$ is of genus at least one, and that $\Sigma$ is an infinite set.
Moreover, for any closed normal subgroup $\mathcal{H}$ of $\Pi$ such that $\mathcal{H}\cap\Delta$ is open in $\Delta$, let $X_{\mathcal{H}}$ be the covering of $X$ corresponding to $\mathcal{H}$, $k_{\mathcal{H}}$ the Galois extension of $k$ corresponding to the image of $\mathcal{H}$ in $G_k$, $G_{k_{\mathcal{H}}}$ the image of $\mathcal{H}$ in $G_k$, and $I_{X_{\mathcal{H}}}$ the set of irreducible component of $X_{\mathcal{H}}$.
Set $I_X:=I_{X_\Pi}$.

\begin{Prop}\label{et}
For any open normal subgroup $\mathcal{H}$ of $\Pi$, the set $I_{X_{\mathcal{H}}}$, together with the natural action of $\Pi/\mathcal{H}$ on $I_{X_{\mathcal{H}}}$, is group-theoretically recovered from $\Delta\subset\Pi$ and $\mathcal{H}$.
Moreover, the natural map $I_{X_\mathcal{H}}\to I_X$ is also recovered.
\end{Prop}

\prf

By Proposition \ref{irred}, the natural map $I_{X_{\mathcal{H}\cap\Delta}}\to I_{X_{\Delta}}$ and the action of $\Delta/(\mathcal{H}\cap\Delta)$ on $I_{X_{\mathcal{H}\cap\Delta}}$ are group-theoretically recovered from $\mathcal{H}\cap\Delta\subset\Delta$.
The assertions follow immediately from the fact that $G_k$ (resp. $G_{k_{\mathcal{H}}}$) naturally acts on $I_{X_{\Delta}}$ (resp. $I_{X_{\mathcal{H}\cap\Delta}}$), and $I_X$ (resp. $I_{X_{\mathcal{H}}}$) is recovered as the set of $G_k$-orbits (resp. $G_{k_{\mathcal{H}}}$-orbits) of $I_{X_{\Delta}}$ (resp. $I_{X_{\mathcal{H}\cap\Delta}}$).
\qed

\

\begin{Rem}\label{decomp}
The natural maps $I_{X_{\mathcal{H}}}\to I_X$ for open normal subgroups $\mathcal{H}$ of $\Pi$ determine a projective limit:
\[\mathfrak{I}:=\varprojlim_{\mathcal{H}\subset\Pi}I_{X_{\mathcal{H}}},\]
where $\mathcal{H}$ runs through the set of open normal subgroups of $\Pi$.
Note that $\Pi$ acts naturally on $\mathfrak{I}$.
For $J=\{J_{\mathcal{H}}\}_{\mathcal{H}\subset\Pi}\in\mathfrak{I}$, let be $\Pi_J$ the stabilizer of $J$ with respect to the action of $\Pi$.
The $\Pi$-conjugacy class of $\Pi_J$ only depends on the irreducible component $J_0:=J_{\Pi}$ of $X$.
Note that the conjugacy class of $\Pi_J$ coincides with that of the subgroup of $\Pi$ determined by the natural morphism $J_0^\nu\to X$, where $J_0^\nu$ is the normalization of $J_0$.
In particular, if $J_0$ is geometrically irreducible over $k$ (or, equivalently, $\Pi_J$ surjects onto $G_k$),
$\Pi_J$ is isomorphic to the (geometrically pro-$\Sigma$) tame fundamental group of $J_0^\nu$.

Proposition \ref{et} says that, we may recover the set $\mathfrak{I}$, together with the natural action of $\Pi$ on $\mathfrak{I}$, is group-theoretically recovered from $\Delta\subset\Pi$.
\end{Rem}

\begin{Def}
For $J=\{J_{\mathcal{H}}\}_{\mathcal{H}\subset\Pi}\in\mathfrak{I}$ and $J_0=J_{\Pi}$, the stabilizer $\Pi_J\subset\Pi$ of $J$ is said to be the {\it{decomposition group of}} $J$.
Moreover, the $\Pi$-conjugacy class of $\Pi_J$ (which depends only on $J_0$) is said to be (the conjugacy class of) the {\it{decomposition group of}} $J_0$.
We shall also denote the decomposition group of $J_0$ (determined up to $\Pi$-conjugacy) by $\Pi_{J_0}$.
\end{Def}

\

In the following, we suppose that, for $\Box\in\{\circ, \bullet\}$, each irreducible component of the normalization of $\overline{X_\Box}$ is of genus at least two, and that $|\Sigma_\circ|, |\Sigma_\bullet|\geq 2$.

\begin{Rem}\label{sigma}
Suppose that there exists an isomorphism of profinite groups $\phi:\Pi_\circ\stackrel{\sim}{\to}\Pi_\bullet$ which induces an isomorphism $\Delta_\circ\stackrel{\sim}{\to}\Delta_\bullet$.
Then, by Remark \ref{p}, it holds that $\Sigma_\circ^\dagger=\Sigma_\bullet^\dagger$.
Moreover, $p_{k_\circ}\in\Sigma_\circ$ if and only if $p_{k_\bullet}\in\Sigma_\bullet$, and if one (hence both) of these conditions holds, it holds that $p_{k_\circ}=p_{k_\bullet}$.
Therefore, it also holds that $\Sigma_\circ=\Sigma_\bullet$.
\end{Rem}

\

By combining the above results with various known results in anabelian geometry, we obtain the following:

\begin{Thm}\label{projective}
Suppose that, for  $\Box\in\{\circ, \bullet\}$, each irreducible component of the normalization of $\overline{X_\Box}$ is of genus at least two.
Suppose, moreover, that there exists an isomorphism of profinite groups $\phi:\Pi_\circ\stackrel{\sim}{\to}\Pi_\bullet$ which induces an isomorphism $\Delta_\circ\stackrel{\sim}{\to}\Delta_\bullet$.
Set 
$\Sigma:=\Sigma_\circ=\Sigma_\bullet$ and $\Sigma':=\Sigma'_\circ=\Sigma'_\bullet$ (cf. Remark \ref{sigma}).
Suppose, furthermore, that $\Sigma$ is an infinite set.
Then the following hold:
\begin{enumerate}
\item[(i)] Suppose that $k_\circ$ and $k_\bullet$ are finite fields.
Suppose, moreover, that $\Sigma'$ is finite.
Then $X_{\circ}^\nu$ and $X_\bullet^\nu$ are isomorphic as schemes.

\item[(ii)] Suppose that $k_\circ=k_\bullet$, and $k:=k_\circ=k_\bullet$ is a finitely generated transcendental extension field of a finite field.
Suppose, moreover, that $p_{k_\circ}=p_{k_\bullet}\not\in\Sigma$, that $\Sigma'$ is finite, that $\phi$ induces the identity automorphism of $G_{k_\circ}=G_{k_\bullet}$, and that any irreducible component of the normalization of $X_\circ\otimes_{k_\circ}k^\mathrm{s}_\circ$ is not isotrivial (i.e., not defined over any finite field).
Then, after replacing the irreducible components of $X_\circ^\nu$ and $X_\bullet^\nu$ by suitable Frobenius twists of them, $X_\circ^\nu$ and $X_\bullet^\nu$ are isomorphic as $k$-schemes.

\item[(iii)] Suppose that $k_\circ$ and $k_\bullet$ are the perfections of finitely generated transcendental extensions of finite fields.
Suppose, moreover, that $p_{k_\circ}, p_{k_\bullet}\not\in\Sigma$, that $\Sigma'$ is finite, and that any irreducible component of the normalization of $X_\circ\otimes_{k_\circ}k^\mathrm{s}_\circ$ is not isotrivial.
Then $X_\circ^\nu$ and $X_\bullet^\nu$ are isomorphic as schemes.

\item[(iv)] Suppose that one of the following holds:
\begin{enumerate}
\item[(a)] $k_\circ$ and $k_\bullet$ are $p$-adic local fields for some $p\in\mathfrak{Primes}$, and $p\in\Sigma$.

\item[(b)] $k_\circ$ and $k_\bullet$ are finitely generated extensions of $\mathbb{Q}$.
\end{enumerate}
Then $X_{\circ}^\nu$ and $X_\bullet^\nu$ are isomorphic as schemes.

\item[(v)] Suppose that $k_\circ=k_\bullet$, and $k:=k_\circ=k_\bullet$ is a generalized sub-$p$-adic field (for some $p\in\mathfrak{Primes}$).
Suppose, moreover, that $p\in\Sigma$, and that $\phi$ induces the identity automorphism of $G_{k_\circ}=G_{k_\bullet}$.
Then $X_{\circ}^\nu$ and $X_\bullet^\nu$ are isomorphic as $k$-schemes.
\end{enumerate}
\end{Thm}

\prf

(i) follows immediately from Proposition \ref{et}, Remark \ref{decomp} and \cite[Theorem C]{ST2}.
(ii) and (iii) follow immediately from Proposition \ref{et}, Remark \ref{decomp} and \cite[Theorems A and B]{Ts} (see also \cite{St1}, \cite{St2}) (note that $p_\circ=p_\bullet$ and $p:=p_\circ=p_\bullet$ is characterized as the unique prime number $l$ such that $\mathrm{cd}_l\,G_\circ=\mathrm{cd}_l\,G_\bullet=1$ (cf. \cite[Proposition 6.5.10 and Theorem 6.5.14]{NSW})).
(iv) follows immediately from Proposition \ref{et}, Remark \ref{decomp}, \cite[Theorem 1.3]{P2},
\cite[Theorem A]{pro-p}, \cite[Corollary 1.3.5]{abs} and \cite[Theorem D]{MoT}.
(v) follows immediately from Proposition \ref{et}, Remark \ref{decomp} and \cite[Theorem 4.12]{sur}.
\qed

\begin{Rem}
Suppose that one of the following conditions holds:

\begin{enumerate}
\item[(i)] $k_\circ$ and $k_\bullet$ are (the perfections of) finitely generated transcendental extensions of finite fields.

\item[(ii)] $k_\circ$ and $k_\bullet$ are finitely generated extensions of $\mathbb{Q}$.
\end{enumerate}
Then any isomorphism $\phi:\Pi_\circ\stackrel{\sim}{\to}\Pi_\bullet$ induces an isomorphism $\Delta_\circ\stackrel{\sim}{\to}\Delta_\bullet$.
Indeed, if (i) or (ii) holds, by a similar argument to the proof of \cite[Lemma 1.1.4 (i)]{abs}, this follows from the fact that $G_{k_\circ}$ and $G_{k_\bullet}$ are very elastic (cf. \cite[Theorem 13.4.2]{FJ} and \cite[Proposition 16.11.6]{FJ}).
The author at the time of writing does not know whether similar results hold over finite fields and $p$-adic local fields or not.
\end{Rem}

\begin{Rem}\label{nonisom}
We can only recover the isomorphism class of $X^\nu$, and cannot recover that of $X$ itself from $\Delta\subset\Pi$.
For example, let $X_\bullet$ be a geometrically irreducible singular curve over a field $k$, $\nu:X_\circ\to X_\bullet$ the normalization of $X_\bullet$ and $\overline{\nu}:=\nu\otimes_k k^\mathrm{s}$.
Suppose that $|\overline{\nu}^{-1}(x)|=1$ for any singular point $x\in \overline{X_\bullet}$.
Then $\nu$ induces an isomorphism of tame fundamental groups $\Pi_\circ\stackrel{\sim}{\to}\Pi_\bullet$, although $X_\circ$ and $X_\bullet$ are not isomorphic.
\end{Rem}

\

\

\section{Tame fundamental groups of non-projective singular curves over finite fields}\label{FF}

For $\Box\in\{\circ, \bullet, \emptyset\}$, let $k_\Box$ be a field of characteristic $p_{k_\Box}\,(\in\{0\}\cup\mathfrak{Primes})$.
Moreover, we shall write:

\begin{itemize}

\item $k_\Box^\mathrm{s}$ for a separable closure of $k_\Box$;

\item $G_{k_\Box}$ for the Galois group $\mathrm{Gal}(k^\mathrm{s}_\Box/k_\Box)$;

\item $X_\Box$ for a projective curve over $k_\Box$ (possibly singular or reducible);

\item $n_\Box$ for the number of irreducible components of $\overline{X_\Box}:=X_\Box\otimes_{k_\Box}k^\mathrm{s}_\Box$;

\item $\nu_\Box: X_\Box^\nu\to X_\Box$ and $\overline{\nu_\Box}:\overline{X_\Box}^\nu\to\overline{X_\Box}$ for the normalization maps;

\item $\overline{X_{\Box, 1}}, \cdots, \overline{X_{\Box, n_\Box}}\subset\overline{X_\Box}$ for the (distinct) irreducible components of $\overline{X_\Box}$;

\item $\overline{X_{\Box, i}}^\nu:=\overline{\nu_\Box}^{-1}(\overline{X_{\Box, i}})\, (1\leq i\leq n_\Box)$;

\item $g_{\Box, i}$ for the genus of $\overline{X_{\Box, i}}^\nu\, (1\leq i\leq n_\Box)$;

\item $\displaystyle g_\Box:=\sum_{i=1}^{n_\Box}g_{\Box, i}$;

\item $D_\Box\subset X_\Box$ for a closed subscheme of $X_\Box$ which is finite \'{e}tale of degree $r_\Box$ over $k_\Box$, and contained in the smooth locus of $X_\Box$;

\item $r_{\Box, i}$ for the degree of $(D_\Box\otimes_{k_\Box}k^\mathrm{s}_\Box)\cap \overline{X_{\Box, i}}$ over $k^\mathrm{s}_\Box$ ($1\leq i\leq n_\Box$);

\item $U_\Box:=X_\Box\setminus D_\Box$;

\item $\overline{U_\Box}:=U_\Box\otimes_{k_\Box}k^\mathrm{s}_\Box$;

\item $\overline{U_{\Box, i}}:=\overline{X_{\Box, i}}\cap \overline{U_\Box}\, (1\leq i\leq n_\Box)$;

\item $\overline{U_{\Box, i}}^\nu:=\overline{\nu_\Box}^{-1}(\overline{U_{\Box, i}})\, (1\leq i\leq n_\Box)$;

\item $n^{>0}_{\Box, \mathrm{p}}:=|\{i\in\mathbb{Z}\,|\,1\leq i\leq n_\Box, g_{\Box, i}>0, r_{\Box, i}=0\}|$;

\item $n_{\Box, \mathrm{a}}:=|\{i\in\mathbb{Z}\,|\,1\leq i\leq n_\Box, r_{\Box, i}>0\}|$;

\item $n^{>0}_{\Box, \mathrm{a}}:=|\{i\in\mathbb{Z}\,|\,1\leq i\leq n_\Box, g_{\Box, i}>0, r_{\Box, i}>0\}|$;

\item $n_{\Box, >2}^{0}:=|\{i\in\mathbb{Z}\,|\,1\leq i\leq n_\Box, g_{\Box, i}=0, r_{\Box, i}>2\}|$;

\item $n_{\Box, 2}^{0}:=|\{i\in\mathbb{Z}\,|\,1\leq i\leq n_\Box, g_{\Box, i}=0, r_{\Box, i}=2\}|$;

\item $n_{\Box, 1}^{0}:=|\{i\in\mathbb{Z}\,|\,1\leq i\leq n_\Box, g_{\Box, i}=0, r_{\Box, i}=1\}|$;

\item $n_{\Box, 0}^{0}:=|\{i\in\mathbb{Z}\,|\,1\leq i\leq n_\Box, g_{\Box, i}=0, r_{\Box, i}=0\}|$;

\item $\Sigma_\Box\subset\mathfrak{Primes}$ for a nonempty set of prime numbers;

\item $\Sigma^\dagger_\Box:=\Sigma_\Box\setminus\{p_{k_\Box}\}$;

\item $\Sigma'_\Box:=\mathfrak{Primes}\setminus\Sigma_\Box$;

\item $\Delta_{X_\Box}:=\pi_1^\mathrm{t}(\overline{X_\Box})^{\Sigma_\Box}$, $\Delta_{X_{\Box, i}}:=\pi_1^\mathrm{t}(\overline{X_{\Box, i}})^{\Sigma_\Box}$, $\Delta_{X^\nu_{\Box, i}}:=\pi_1^\mathrm{t}(\overline{X_{\Box, i}}^\nu)^{\Sigma_\Box}$, $\Delta_{U_\Box}:=\pi_1^\mathrm{t}(\overline{U_\Box})^{\Sigma_\Box}$, $\Delta_{U_{\Box, i}}:=\pi_1^\mathrm{t}(\overline{U_{\Box, i}})^{\Sigma_\Box}$, $\Delta_{U^\nu_{\Box, i}}:=\pi_1^\mathrm{t}(\overline{U_{\Box, i}}^\nu)^{\Sigma_\Box}$;

\item $\Pi_{X_\Box}:=\pi_1^\mathrm{t}(X_\Box)/\mathrm{Ker}(\pi_1^\mathrm{t}(\overline{X_\Box})\twoheadrightarrow\Delta_{X_\Box})$, $\Pi_{U_\Box}:=\pi_1^\mathrm{t}(U_\Box)/\mathrm{Ker}(\pi_1^\mathrm{t}(\overline{U_\Box})\twoheadrightarrow\Delta_{U_\Box})$.

\end{itemize}

In the following, for $\Box\in\{\circ, \bullet, \emptyset\}$, we suppose that $g_\Box+r_\Box-n_{\Box, \mathrm{a}}>0$.
Moreover, suppose that we are given not only $\Pi_{U_\Box}$ but also the subgroup $\Delta_{U_\Box}$ as the input datum (i.e., we reconstruct various invariants of $U_\Box$ group-theoretically from $\Delta_{U_\Box}\subset\Pi_{U_\Box}$).
In other words, we are given the following exact sequence as the input data:
\[\xymatrix{1 \ar[r]  & \Delta_{U_\Box} \ar[r] & \Pi_{U_\Box} \ar[r]^{\mathrm{pr}_{U_\Box}} & G_{k_\Box} \ar[r] & 1.}\]
($\Pi_{U_\Box}$ surjects onto $G_{k_\Box}$ since $X_\Box$ (hence also $U_\Box$) is geometrically connected.)

Note that, by a similar argument to \cite[\S 4.2]{D}, we have the following isomorphism of profinite groups:
\begin{align}\label{tisom}
\Delta_{U_\Box}\simeq\left(\Star_{i=1}^{n_\Box}\Delta_{U_{\Box, i}^\nu}\right)\ast F_{\delta_\Box},
\end{align}
where $\delta_\Box:=\displaystyle 1-n_\Box+\sum_{x\in \overline{X}}\left(|(\overline{\nu_\Box})^{-1}(x)|-1\right)$ and $F_{\delta_\Box}$ is a free pro-$\Sigma_\Box$ group of rank $\delta_\Box$.

Set $\Delta_{\Box, i}:=\Delta_{U_{\Box, i}^\nu}$ for $1\leq i\leq n_\Box$ and $\Delta_\Box:=\displaystyle\left(\Star_{i=1}^{n_\Box} \Delta_{\Box, i}\right)\ast F_{\delta_\Box}$, and fix an isomorphism $\iota_\Box:\Delta_{U_\Box}\stackrel{\sim}{\to}\Delta_\Box$ (cf. (\ref{tisom})).

\begin{Def}
For any open normal subgroup $\mathcal{H}$ of $\Pi_U$, let $I_{U_\mathcal{H}}$ be the set of irreducible components of $U_\mathcal{H}$.
As in Remark  \ref{decomp}, for open normal subgroups $\mathcal{H}'\subset\mathcal{H}$ of $\Pi_U$, there exists a natural map $I_{U_{\mathcal{H}'}}\to I_{U_\mathcal{H}}$, and these maps determine a projective limit:
\[\mathfrak{I}:=\lim_{\mathcal{H}\subset\Pi_U}I_{U_\mathcal{H}},\]
where $\mathcal{H}$ runs through the set of open normal subgroups of $\Pi_U$.
$\Pi_U$ acts naturally on $\mathfrak{I}$.
Suppose that $J=\{J_{\mathcal{H}}\}_{\mathcal{H}\subset\Pi_U}\in\mathfrak{I}$ and $J_0:=J_{\Pi_U}$.
Then the stabilizer $\Pi_J\subset\Pi_U$ of $J$ is said to be the {\it{decomposition group of}} $J$.
Moreover, the $\Pi_U$-conjugacy class of $\Pi_J$ (which depends only on $J_0$) is said to be the (conjugacy class of) {\it{decomposition group of}} $J_0$.
We shall also denote the decomposition group of $J_0$ (determined up to $\Pi_U$-conjugacy) by $\Pi_{J_0}$.
\end{Def}

\begin{Rem}\label{afdecomp}
Let us take an element $J_0$ of $I_U:=I_{U_{\Pi_U}}$.
As in Remark \ref{decomp}, if $J_0$ is geometrically irreducible over $k$ (or, equivalently, $\Pi_{J_0}$ surjects onto $G_k$),
the subgroup $\Pi_{J_0}$ of $\Pi_U$ (determined up to conjugacy) is isomorphic to the (geometrically pro-$\Sigma$) tame fundamental group of $J_0^\nu$, where $J_0^\nu$ is the inverse image of $J_0$ by $\nu:X^\nu\to X$.
\end{Rem}

\

Until the end of this section, suppose that $k$ is a finite field of characteristic $p_k$.
Set $q_k:=|k|$, and let $\varphi_k\in G_k$ be the $q_k$-th power Frobenius element.
Suppose, moreover, that $\Sigma'$ is a finite set.

\begin{Rem}\label{pfin}
By a similar argument to Remark \ref{p}, we may determine $\Sigma^\dagger$ group-theoretically from $\Delta_U$.
Moreover, if $\overline{U_i}^\nu$ is a hyperbolic curve for some $1\leq i\leq n$, or if $\delta\neq 0$, we may also determine group-theoretically whether $\Sigma$ contains $p_k$ or not, and $p_k$ itself if $p_k\in\Sigma$.
If $p_k\in\Sigma$, $p_k$ is characterized as the unique prime number $l$ satisfying the following condition (cf. the assumption that $g+r-n_\mathrm{a}>0$):
\begin{quote}
There exists an open subgroup $\mathcal{H}$ of $\Delta_U$ such that $\mathcal{H}^{l\text{-ab}}$ is not trivial, and that 
$\mathrm{rank}_{\mathbb{Z}_l}\mathcal{H}^{l\text{-ab}}<\mathrm{rank}_{\mathbb{Z}_{l'}}\mathcal{H}^{l'\text{-ab}}$ for any $l'\in\Sigma^\dagger$
\end{quote}
\end{Rem}

\

By a similar argument to the proof of \cite[Proposition 3.4]{Tama1} and \cite[Proposition 2.5]{ST2}, we may recover $p_k$ (even if $p_k\not\in\Sigma$ (cf. Remark \ref{pfin})), $q_k$ and the ($q_k$-th power) Frobenius element of $G_k$:

\begin{Prop}\label{Frob}
$p_k$, $q_k$ and the Frobenius element $\varphi_k\in G_k=\Pi_U/\Delta_U$ are group-theoretically recovered from $\Delta_U\subset\Pi_U$.

\end{Prop}

\prf

Consider the natural character
\[\rho^\mathrm{det}:G_k\to\mathrm{Aut}\left(\bigwedge_{\hat{\mathbb{Z}}^{\Sigma^\dagger}}^{\mathrm{max}}\left(\Delta_U^{\mathrm{ab}}\right)^{\Sigma^\dagger}\right)=\left(\hat{\mathbb{Z}}^{\Sigma^\dagger}\right)^\times.\]
Note that $\rho^{\mathrm{det}}$ may be recovered from $\Delta_U\subset\Pi_U$ group-theoretically (cf. Remark \ref{pfin}).
Let $N_0$ be a positive integer such that each irreducible component of $X_{k_{N_0}}:=X\otimes_k k_{N_0}$ is geometrically irreducible, where $k_{N_0}$ is the extension of $k$ of degree $N_0$.
Then we have $(\rho^{\mathrm{det}})^2|_{G_{k_{N_0}}}=(\chi_{\Sigma^\dagger, N_0})^{2(g+r-n_{\mathrm{a}})}$, where $G_{k_{N_0}}\subset G_k$ is the absolute Galois group of $k_{N_0}$ and $\chi_{\Sigma^\dagger, N_0}:G_{k_{N_0}}\to\left(\hat{\mathbb{Z}}^{\Sigma^\dagger}\right)^\times$ is the $\Sigma^\dagger$-part of the cyclotomic character (cf. \cite[Proposition 3.4]{Tama1} and \cite[Theorem 1]{AP2}).
For each positive integer $N$, let $w_N$ be the cardinality of the coinvariant quotient $\hat{\mathbb{Z}}^{\Sigma^\dagger}_{(\rho^\mathrm{det})^2(G_{k_N})}$, where $G_{k_N}\subset G_k$ is the absolute Galois group of the extension of $k$ of degree $N$.
Moreover, for each positive integer $m$, set:
\begin{align*}
\mathcal{M}_m&:=\{M\in\mathbb{R}\,|\, M>0, \exists C>0, \forall N\in\mathbb{N}, w_{mN}\leq CM^{mN}\}, \\
M_{m}&:=\inf\mathcal{M}_m.
\end{align*}
Then, by the arguments of the proof of \cite[Proposition 2.5]{ST2}, if $N_0$ divides $m$, then $M_{m}=q_k^{2m(g+r-n_{\mathrm{a}})}$.
(Note that $\Sigma'$ is finite and, in particular, $\Sigma$ is $k$-large (cf. \cite[Definition 2.1]{ST2}).)
In particular, we have:
\[\lim_{m\to\infty}(M_{m!})^{\frac{1}{m!}}=q_k^{2(g+r-n_{\mathrm{a}})}.\]
Therefore, as in the proof of \cite[Proposition 2.5]{ST2}, $p_k$, $q_k$ and $\varphi_k$ are group-theoretically recovered from $\Delta_U\subset\Pi_U$.
\qed

\

For $1\leq i\leq n$ and some assertion(s) $P_i$ concerning $i$, $\displaystyle\sum_{P_i}$ means to take sum for $i$ satisfying $P_i$.
Moreover, for a closed normal subgroup $\mathcal{H}\subset\Pi_U$ such that $\mathcal{H}\cap\Delta_U$ is open in $\Delta_U$, we denote the covering of $U$ corresponding to $\mathcal{H}$ by $U_\mathcal{H}$, and the extension of $k$ corresponding to the image of $\mathcal{H}$ in $G_k=\Pi_U/\Delta_U$ by $k_{\mathcal{H}}$.
We denote invariants of $U_{\mathcal{H}}/k_{\mathcal{H}}$ by the corresponding invariants of $U/k$ with subscript $\mathcal{H}$ (e.g., $g_{\mathcal{H}}$, $n_{\mathcal{H}}$, $(n_\mathrm{p}^{>0})_\mathcal{H}$ etc.)

\begin{Prop}\label{inv}
$g$, $n^{>0}_\mathrm{p}$, $n^{>0}_\mathrm{a}$, $n^0_{>2}$, $n^0_{2}$, $r-n^0_1$ and $\delta$ are group-theoretically recovered from $\Delta_U\subset\Pi_U$.
\end{Prop}

\prf

If $p_k\in\Sigma$, by replacing $\Delta_U$ by its maximal pro-$\Sigma^\dagger$ quotient, we assume that $p_k\not\in\Sigma$.

Fix $l\in\Sigma$ (hence $l\neq p$), and let $P$ be the characteristic polynomial of the $q_k$-th power Frobenius element $\varphi_k\in G_k$ (which is group-theoretically determined from $\Delta_U\subset\Pi_U$ (cf. Proposition \ref{Frob})) on the free $\mathbb{Z}_l$-module $\Delta_U^{l\text{-ab}}$.
Then, $2g$ (resp. $r-n_\mathrm{a}$, resp. $\delta$) is group-theoretically recovered as the number of roots of $P$ (with multiplicity) whose absolute values are $\sqrt{q_k}$ (resp. $q_k$, resp. 1).
Moreover, similarly to the proof of Proposition \ref{delta}, we obtain $n^{>0}_\mathrm{p}=\mathrm{dim}_{\mathbb{Z}/l\mathbb{Z}}H^2(\Delta_U^l, \mathbb{Z}/l\mathbb{Z})$.

In the following, we reconstruct $n_\mathrm{a}^{>0}$, $n^0_{>2}$, $n^0_2$ and $r-n_1^0$ group-theoretically.
Set:
\begin{align*}
\mathfrak{H}&:=\{\mathcal{H}\subset\Pi_U:\text{ an open normal subgroup}\,|\,(\Delta_U:\Delta_U\cap\mathcal{H})=l\}, \\
M_1&:=\min\{2g_\mathcal{H}\,|\,\mathcal{H}\in\mathfrak{H}\}, \\
\mathfrak{H}_1&:=\{\mathcal{H}\in\mathfrak{H}\,|\,2g_\mathcal{H}=M_1\}.
\end{align*}
Note that these are group-theoretically recovered from $\Delta_U\subset\Pi_U$.
Here, for $\mathcal{H}\in\mathfrak{H}$, $\mathcal{H}$ belongs to $\mathfrak{H}_1$ if and only if, for $i$ satisfying $g_i>0$ (resp. for $i$ satisfying $g_i=0$ and $r_i>1$), $f_i$ corresponds (via $\iota$) to a nontrivial finite \'{e}tale covering of $\overline{U_i}^\nu$ unramified over $\overline{X_i}^\nu\setminus \overline{U_i}^\nu$ (resp. a finite \'{e}tale covering of $\overline{U_i}^\nu$ which is trivial or ramified at only two points of $\overline{X_i}^\nu\setminus\overline{U_i}^\nu$), where 
\[f=(f_i)_{i=1}^n+f_\delta\in\displaystyle\left(\bigoplus_{i=1}^n H^1(\Delta_i, \mathbb{Z}/l\mathbb{Z})\right)\oplus H^1(F_\delta, \mathbb{Z}/l\mathbb{Z})=H^1(\Delta, \mathbb{Z}/l\mathbb{Z})\] is the element corresponding to the natural surjection $\Delta\twoheadrightarrow \Delta/\iota(\Delta_U\cap\mathcal{H})$.
Thus, by the Hurwitz formula, for any $\mathcal{H}\in\mathfrak{H}_1$, we have:
\begin{align*}
(M_1=)\,2g_\mathcal{H}&=\sum_{g_i>0}(l(2g_i-2)+2)=2lg-2(l-1)(n^{>0}_\mathrm{p}+n_{\mathrm{a}}^{>0}).
\end{align*}
Therefore, $n_{\mathrm{a}}^{>0}$ is recovered group-theoretically.
Moreover,
\begin{align*}
\min\{r_\mathcal{H}-(n_\mathrm{a})_\mathcal{H}\,|\,\mathcal{H}\in\mathfrak{H}_1\}&=\sum_{g_i>0, r_i>0}(lr_i-1)+\sum_{g_i=0, r_i>1}(l(r_i-2)+2-1)+\sum_{(g_i, r_i)=(0, 1)}(lr_i-l) \\
&=lr-n^{>0}_\mathrm{a}-(2l-1)(n^0_{>2}+n^0_2)-ln^0_1.
\end{align*}
(Note that, for $\mathcal{H}\in\mathfrak{H}_1$, $r_\mathcal{H}-(n_\mathrm{a})_\mathcal{H}$ attains the minimum if and only if, for $i$ satisfying $g_i=0$ and $r_i>1$, $f_i$ corresponds (via $\iota$) to a finite \'{e}tale covering of $\overline{U_i}^\nu$ ramified at only two points of $\overline{X_i}^\nu\setminus \overline{U_i}^\nu$.)
So, $l(r-n^0_1)-(2l-1)(n^0_{>2}+n^0_2)$ is recovered group-theoretically.

On the other hand, if $n^{>0}_\mathrm{p}=n^{>0}_\mathrm{a}=\delta=0$ (note that this is a group-theoretic condition as above), then $n^0_{>2}+n^0_2>0$ (cf. the assumption that $g+r-n_\mathrm{a}>0$).
In this case, there exists a positive integer $1\leq j\leq n$ such that $g_j=0$, $r_j>1$, and we have
\begin{align*}
\max\{r_\mathcal{H}-(n_\mathrm{a})_\mathcal{H}\,|\,\mathcal{H}\in\mathfrak{H}_1\}&=\sum_{g_i=0, r_i>1, i\neq j}(lr_i-l)+(l(r_j-2)+2-1)+\sum_{(g_i, r_i)=(0, 1)}(lr_i-l) \\
&=lr-l(n^0_{>2}+n^0_2)-l+1-ln^0_1,
\end{align*}
and otherwise (i.e., if at least one of $n^{>0}_\mathrm{p}$, $n^{>0}_\mathrm{a}$ and $\delta$ is nonzero),
\begin{align*}
\max\{r_\mathcal{H}-(n_\mathrm{a})_\mathcal{H}\,|\,\mathcal{H}\in\mathfrak{H}_1\}&=\sum_{g_i>0, r_i>0}(lr_i-1)+\sum_{g_i=0, r_i>0}(lr_i-l)+\sum_{(g_i, r_i)=(0, 1)}(lr_i-l) \\
&=lr-n^{>0}_\mathrm{a}-l(n^0_{>2}+n^0_2+n^0_1).
\end{align*}
Note that, in the case where $n^{>0}_\mathrm{p}=n^{>0}_\mathrm{a}=\delta=0$ (resp. at least one of $n^{>0}_\mathrm{p}$, $n^{>0}_\mathrm{a}$ and $\delta$ is nonzero), for $\mathcal{H}\in\mathfrak{H}_1$, $r_\mathcal{H}-(n_\mathrm{a})_\mathcal{H}$ attains the maximum if and only if, there exists $j$ such that (necessarily $g_j=0$ and) $r_j>1$ and that $f_j$ corresponds (via $\iota$) to a finite \'{e}tale covering of $\overline{U_j}^\nu$ ramified at only two points of $
\overline{X_j}^\nu\setminus \overline{U_j}^\nu$, and $f_i$ is trivial for $i$ satisfying (necessarily $g_i=0$ and) $r_i>1$ and $i\neq j$ (resp. $f_i$ is trivial for $i$ satisfying $g_i=0$ and $r_i>1$).
In both cases, $l(r-n^0_1)-l(n^0_{>2}+n^0_2)$ is recovered group-theoretically.
Hence also $r-n^0_1$ and $n^0_{>2}+n^0_2$ are recovered group-theoretically.

Finally, we shall determine $n^0_{>2}$ and $n^0_2$.
If $n^0_{>2}+n^0_2=0$, then we have $n^0_{>2}=n^0_2=0$.
So, in the following, we assume that $n^0_{>2}+n^0_2>0$.
Moreover, suppose that $l\neq 2$ (note that $\Sigma$ is an infinite set).
Set:
\begin{align*}
\mathfrak{H}_2&:=\{\mathcal{H}\in\mathfrak{H}\,|\,(n^{>0}_\mathrm{p})_\mathcal{H}=ln^{>0}_\mathrm{p}, (n^{>0}_\mathrm{a})_{\mathcal{H}}=\max_{\mathcal{H}'\in\mathfrak{H}}(n^{>0}_\mathrm{a})_{\mathcal{H}'}\}, \\
M_2&:=\max\{2g_\mathcal{H}\,|\,\mathcal{H}\in\mathfrak{H}_2\}, \\
\mathfrak{H}_3&:=\{\mathcal{H}\in\mathfrak{H}_2\,|\,2g_\mathcal{H}=M_2\}.
\end{align*}
(Note that, $\mathcal{H}\in\mathfrak{H}$ satisfies $\displaystyle (n^{>0}_\mathrm{a})_{\mathcal{H}}=\max_{\mathcal{H}'\in\mathfrak{H}}(n^{>0}_\mathrm{a})_{\mathcal{H}'}$ if and only if $f_i$ is trivial for $i$ satisfying $g_i>0$ and $r_i>0$, and corresponds to a finite \'{e}tale covering of $\overline{U_i}^\nu$ ramified at three or more points of $\overline{X_i}^\nu\setminus\overline{U_i}^\nu$ for $i$ satisfying $g_i=0$ and $r_i>2$.
Such covering exists since we have assumed that $n^0_{>2}+n^0_2>0$ and $l\neq 2$.)
Then, by the Hurwitz formula, we have:
\[M_2=2lg+\sum_{g_i=0, r_i>2}(l-1)(r_i-2)=2lg+(l-1)\sum_{g_i=0, r_i>2}r_i-2(l-1)n^0_{>2}.\]
Note that $M_2=2lg$ if and only if $n^0_{>2}=0$ (and hence $n^0_2=n^0_{>2}+n^0_2$).
Therefore, we assume that $M_2>2lg$ (i.e., $n^0_{>2}>0$).
Then, we have:
\begin{align*}
M_3&:=\max\{r_\mathcal{H}-(n_\mathrm{a})_\mathcal{H}\,|\,\mathcal{H}\in\mathfrak{H}_3\} \\
&=\sum_{g_i>0, r_i>0}(lr_i-l)+\sum_{g_i=0, r_i>2}(r_i-1)+\sum_{g_i=0, r_i=2}(lr_i-l)+\sum_{g_i=0, r_i=1}(lr_i-l) \\
&=lr-(l-1)\sum_{g_i=0, r_i>2}r_i-ln^{>0}_\mathrm{a}-n^0_{>2}-ln^0_2-ln^0_1
\end{align*}
Since $M_2$, $M_3$, $2g$, $n^{>0}_\mathrm{a}$, $r-n^0_1$ and $n^0_{>2}+n^0_2$ have been already recovered group-theoretically, $n^0_{>2}$ and $n^0_2$ are also recovered, as desired.
\qed

\begin{Cor}\label{surj}
The natural surjections $\Delta_U\twoheadrightarrow\Delta_X$ and $\Delta\twoheadrightarrow F_\delta$ are group-theoretically recovered from $\Delta_U\subset\Pi_U$ (and $\iota$).
\end{Cor}

\prf

Let $\mathfrak{H}$ be the set of open normal subgroups of $\Delta_U$ consisting of elements $\mathcal{H}\subset\Delta_U$ such that $r_\mathcal{H}=(\Delta_U:\mathcal{H})r$.
Note that, although $r$ itself has not been group-theoretically recovered yet, for any open normal subgroup $\mathcal{H}\subset\Delta_U$, we have necessarily $(n^0_1)_\mathcal{H}=(\Delta_U:\mathcal{H})n^0_1$.
Therefore, $r_\mathcal{H}=(\Delta_U:\mathcal{H})r$ if and only if $r_\mathcal{H}-(n^0_1)_\mathcal{H}=(\Delta_U:\mathcal{H})(r-n^0_1)$, and the latter is a group-theoretic condition by Proposition \ref{inv}.

It is clear that
\[\Delta_X=\Delta_U/\bigcap_{\mathcal{H}\in\mathfrak{H}}\mathcal{H}.\]
Moreover, note that the surjection $\Delta\twoheadrightarrow F_\delta$ necessarily factors through $\overline{\iota}(\Delta_X)$, where $\overline{\iota}$ is the isomorphism $\overline{\iota}:\Delta_X\stackrel{\sim}{\to}\Delta/\iota(\mathrm{Ker}(\Delta_U\twoheadrightarrow\Delta_X))$ induced by $\iota$.
The kernel of the surjection $\overline{\iota}(\Delta_X)\twoheadrightarrow F_\delta$ is recovered as in the proof of Proposition \ref{delta}.
\qed

\

In the following, we assume that $(g_i, r_i)\not\in\{(0, 0), (0, 1)\}$ for $1\leq i\leq n$ (i.e., $n^0_1=n^0_0=0$).
Then we have a similar result to Proposition \ref{et}:

\begin{Prop}\label{irredfin}
For any open normal subgroup $\mathcal{H}$ of $\Pi_U$, the set $I_{U_{\mathcal{H}}}$, together with the natural action of $\Pi_U/\mathcal{H}$ on $I_{U_{\mathcal{H}}}$, is group-theoretically recovered from $\Delta_U\subset\Pi_U$.
In particular, $\mathfrak{I}=\displaystyle\varprojlim_{\mathcal{H}\subset\Pi_U}I_{U_{\mathcal{H}}}$, together with the natural action of $\Pi_U$ on $\mathfrak{I}$, is recovered group-theoretically.
Moreover, the natural map $I_{U_\mathcal{H}}\to I_U$ is also recovered.
\end{Prop}

\prf

By similar arguments to the proofs of Propositions \ref{irred} and \ref{et}, for any open normal subgroup $\mathcal{H}$ of $\Pi_U$, the natural map $I_{U_{\mathcal{H}}}\to I_U$ is recovered group-theoretically from $\Delta_U\subset\Pi_U$.
(Although Proposition \ref{irred} treats only projective curves such that the genus of each irreducible component of the normalizations of curves is at least one, a similar argument works since we have assumed that $(g_i, r_i)\not\in\{(0, 0), (0, 1)\}$ for $1\leq i\leq n$ (and hence $\Delta_{U_i^\nu}$ is not trivial for $1\leq i\leq n$.)
The point is that we have recovered $n$ and the natural surjection $\Delta\twoheadrightarrow F_\delta$ group-theoretically.)
\qed

\begin{Cor}\label{rat}
$|X(k)|$ and $|U(k)|$ are group-theoretically recovered from $\Delta_U\subset\Pi_U$.
\end{Cor}

\prf

Note that we have already recovered the $q_k$-th power Frobenius element $\varphi_k\in G_k=\Pi_U/\Delta_U$.
Let $I_U'\subset I_U$ be the subset of $I_U$ consisting of irreducible components of $U$ which are geometrically irreducible.
(The geometrically irreducibleness of each $J\in I_U$ is a group-theoretic condition (cf. Remark \ref{afdecomp}).)
Set $n':=|I_U'|$.
For each $J\in I_U'$, let $\Pi_J$ be the decomposition group of $J$ (determined up to conjugacy), $J^\nu$ (resp. $J^\nu_\mathrm{c}$) the inverse image of $J$ (resp. the closure of $J$ in $X$) by $\nu$, $\Delta_{J^\nu}:=\Pi_J\cap \Delta_U$ and $\Delta_{J^\nu_\mathrm{c}}$ the image of $\Delta_{J^\nu}$ by the natural surjection $\Delta_U\twoheadrightarrow \Delta_X$.
Note that, $\Delta_{J^\nu}$ (resp. $\Delta_{J^\nu_\mathrm{c}}$) is isomorphic to the geometric pro-$\Sigma$ tame fundamental group of $J^\nu$ (resp. $J^\nu_\mathrm{c}$), and the outer actions of $G_k$ on $\Delta_{J^\nu}$ and $\Delta_{J^\nu_\mathrm{c}}$ have been already recovered group-theoretically.
Moreover, let $\Delta_\delta^\mathrm{ab}$ be the subgroup of $\Delta_X^\mathrm{ab}$ corresponding to $F_\delta^\mathrm{ab}\subset \Delta^\mathrm{ab}$ via $\iota$.
Then, by the Lefschetz trace formula (see also \cite[Propositions 6, 7 and Theorem 11]{AP2} and \cite[Proposition 1.15]{ST1}), for any $l\in\Sigma^\dagger$,
\begin{align*}
|X(k)|&=\sum_{j=0}^2(-1)^j\mathrm{Tr}(\varphi_k^{-1}\,|\,H^j_{{\mathrm{\'{e}t}}}(\overline{X}, \mathbb{Q}_l)) \\
&=1+n'q_k-\mathrm{Tr}(\varphi_k^{-1}\,|\,H^1_{{\mathrm{\'{e}t}}}(\overline{X}, \mathbb{Q}_l)) \\
&=1+n'q_k-\mathrm{Tr}(\varphi_k^{-1}\,|\,H^1(\Delta_X, \mathbb{Q}_l)), \\
&=1+n'q_k-\mathrm{Tr}(\varphi_k^{-1}\,|\,H^1(\Delta_\delta^\mathrm{ab}, \mathbb{Q}_l))-\sum_{J\in I_U'}\mathrm{Tr}(\varphi_k^{-1}\,|\,H^1(\Delta_{J^{\nu}_\mathrm{c}}, \mathbb{Q}_l))
\end{align*}
where, for a profinite group $\Gamma$, we define
\[H^j(\Gamma, \mathbb{Q}_l):=\left(\varprojlim_{n\to\infty}H^j(\Gamma, \mathbb{Z}/l^n\mathbb{Z})\right)\otimes_{\mathbb{Z}_l}\mathbb{Q}_l.\]
(Note that we have already recovered $n'$, $q_k$ and the surjection $\Delta_U\twoheadrightarrow \Delta_X$ group-theoretically.)
Moreover, by a similar argument to the proof of \cite[Proposition 1.15]{ST1}, for each $J\in I_U'$, the number $|J\cap D(k)|$ is recovered group-theoretically from $\Pi_J$.
Hence,
\[|U(k)|=|X(k)|-\sum_{J\in I_U'}|J\cap D(k)|,\]
is also group-theoretically recovered.
(Note that $D$ is contained in the smooth locus of $X$.)
\qed

\begin{Rem}
In the situation of Corollary \ref{rat}, $|U(k)|$ is also written in the following way:
\begin{align*}
|U(k)|&=|X(k)|-|D(k)| \\
&=|X(k)|-\sum_{J\in I_U'}\sum_{j=0}^2(-1)^j\left(\mathrm{Tr}(\varphi_k^{-1}\,|\,H^j_{\mathrm{\'{e}t}}(\overline{J^\nu_\mathrm{c}}, \mathbb{Q}_l))-\mathrm{Tr}(\varphi_k^{-1}\,|\,H_\mathrm{c}^j(\overline{J^\nu}, \mathbb{Q}_l))\right) \\
&=|X(k)|-\sum_{J\in I_U'}\left(1+q_k-\mathrm{Tr}(\varphi_k^{-1}\,|\,H^1_{\mathrm{\'{e}t}}(\overline{J^\nu_\mathrm{c}}, \mathbb{Q}_l))-\sum_{j=0}^2(-1)^j\mathrm{Tr}(\varphi_k\,|\,H_\mathrm{\'{e}t}^j(\overline{J^\nu}, \mathbb{Q}_l(1)))\right) \\
&=|X(k)|-\sum_{J\in I_U'}\left(1+q_k-\mathrm{Tr}(\varphi_k^{-1}\,|\,H^1_{\mathrm{\'{e}t}}(\overline{J^\nu_\mathrm{c}}, \mathbb{Q}_l))-q_k\sum_{j=0}^2(-1)^j\mathrm{Tr}(\varphi_k\,|\,H_\mathrm{\'{e}t}^j(\overline{J^\nu}, \mathbb{Q}_l))\right) \\
&=|X(k)|-\sum_{J\in I_U'}\left(1+q_k-\mathrm{Tr}(\varphi_k^{-1}\,|\,H^1(\Delta_{J^\nu_\mathrm{c}}, \mathbb{Q}_l))-q_k\sum_{j=0}^2(-1)^j\mathrm{Tr}(\varphi_k\,|\,H^j(\Delta_{J^\nu}, \mathbb{Q}_l))\right) \\
&=1-n'-\mathrm{Tr}(\varphi_k^{-1}\,|\,H^1(\Delta_\delta^\mathrm{ab}, \mathbb{Q}_l))+q_k\sum_{J\in I_U'}\sum_{j=0}^2(-1)^j\mathrm{Tr}(\varphi_k\,|\,H^j(\Delta_{J^\nu}, \mathbb{Q}_l)),
\end{align*}
where $\overline{J^\nu}:=J^\nu\otimes_k k^\mathrm{s}$,  $\overline{J^\nu_\mathrm{c}}:=J^\nu_\mathrm{c}\otimes_k k^\mathrm{s}$.
\end{Rem}

\

By using the above results, we shall consider anabelian phenomena for curves which are not necessarily projective over finite fields.

In the following, we suppose that, for $\Box\in\{\circ, \bullet\}$, $k_\Box$ is a finite field of characteristic $p_{k_\Box}$ such that $|k_\Box|=q_{k_\Box}$.
Moreover, we suppose that each irreducible component of $\overline{U_\Box}^\nu:=(\overline{\nu_\Box})^{-1}(\overline{U_\Box})$ is a curve of type distinct from $(0, 0), (0, 1)$ and $(1, 0)$ (hence $g_\Box+r_\Box-n_{\Box, \mathrm{a}}>0$ automatically holds), and that $\Sigma_\Box$ is an infinite set.

\begin{Rem}\label{sigmafin}
Suppose that there exists an isomorphism of profinite groups $\phi:\Pi_{U_\circ}\stackrel{\sim}{\to}\Pi_{U_\bullet}$ which induces an isomorphism $\Delta_{U_\circ}\stackrel{\sim}{\to}\Delta_{U_\bullet}$.
Then, by Remark \ref{pfin}, it holds that $\Sigma_\circ^\dagger=\Sigma_\bullet^\dagger$.
Moreover, by Proposition \ref{Frob}, we have $p_{k_\circ}=p_{k_\bullet}$ and $q_{k_\circ}=q_{k_\bullet}$ (hence $k_\circ$ and $k_\bullet$ are isomorphic).
However, $\Sigma_\circ=\Sigma_\bullet$ does not necessarily hold.
Indeed, for $\Box\in\{\circ, \bullet\}$, suppose that $U_\Box$ is a smooth curve of type $(0, 2)$ over $k_\Box$, and that $\Sigma_\circ=\mathfrak{Primes}, \Sigma_\bullet=\mathfrak{Primes}\setminus\{p\}$, where $p:=p_{k_\circ}=p_{k_\bullet}$.
Then there exists an isomorphism of profinite groups $\Pi_\circ\stackrel{\sim}{\to}\Pi_\bullet$ inducing an isomorphism $\Delta_\circ\stackrel{\sim}{\to}\Delta_\bullet$.

Of course, if, for some $\Box\in\{\circ, \bullet\}$, $\overline{U_\Box}^\nu$ has an irreducible component which is hyperbolic, or ``$\delta_\Box\neq 0$'' ($\delta_\Box$ is defined as above), then $p_{k_\circ}\in\Sigma_\circ$ if and only if $p_{k_\bullet}\in\Sigma_\bullet$, and hence $\Sigma_\circ=\Sigma_\bullet$.
\end{Rem}

\begin{Thm}\label{mainfin}
Suppose that there exists an isomorphism of profinite groups $\phi:\Pi_{U_\circ}\stackrel{\sim}{\to}\Pi_{U_\bullet}$ which induces an isomorphism $\Delta_{U_\circ}\stackrel{\sim}{\to}\Delta_{U_\bullet}$.
Suppose, moreover, that $\Sigma_\circ'$ and $\Sigma_\bullet'$ are finite sets.
Then $U_\circ^\nu$ and $U_\bullet^\nu$ are isomorphic as schemes.
\end{Thm}

\prf

Immediate from \cite[Theorem D]{ST2}, Remark \ref{afdecomp} and Proposition \ref{irredfin}.
Note that, even if there exists an irreducible component $J_\Box$ of $U_\Box$ for $\Box\in\{\circ, \bullet\}$ which is not geometrically irreducible over $k_\Box$, the decomposition subgroup $\Pi_{J_\Box}\subset\Pi_{U_\Box}$ are recovered (up to conjugacy) group-theoretically.
So, we may recover the minimal finite (Galois) extension $k'_\Box$ of $k_\Box$ such that each irreducible component of $J_\Box\otimes_{k_\Box}k'_\Box$ is geometrically irreducible as the field corresponding to the (open) image of $\Pi_{J_\Box}$ in $G_{k_\Box}$.
$\Pi_{J_\Box}$ is isomorphic to the (geometrically pro-$\Sigma_\Box$) tame fundamental group of the normalization $(J_\Box')^\nu$ of a irreducible component $J_\Box'$ of $J_\Box\otimes_{k_\Box}k'_\Box$, and $J_\Box\otimes_{k_\Box}k'_\Box$ (resp. $J_\Box^\nu\otimes_{k_\Box}k'_\Box$) is isomorphic to the disjoint union of $[k_\Box':k_\Box]$ copies of $J_\Box'$ (resp. $(J_\Box')^\nu$).
Moreover, since the natural morphism $(J'_\Box)^\nu\to J_\Box^\nu$ is finite \'{e}tale of degree one, these are isomorphic as schemes.

Although curves of type $(0, 2)$ (over $k_\Box$, $\Box\in\{\circ, \bullet\}$) are not hyperbolic, these curves are isomorphic (over $k_\Box$) to $\mathbb{P}_{k_\Box}^1\setminus\{0, \infty\}$ or $\mathbb{P}_{k_\Box}^1\setminus\{\alpha_\Box, \overline{\alpha_\Box}\}$, where $\alpha_\Box\in k^\mathrm{s}_\Box$, $[k_\Box(\alpha_\Box):k_\Box]=2$, and $\overline{\alpha_\Box}$ is the conjugate of $\alpha_\Box$ over $k_\Box$.
(Note that the Brauer groups of finite fields are trivial.)
For each irreducible component $J_\Box$ of $U_\Box$ whose normalization is of type $(0, 2)$, let $\Pi_{J_\Box}\subset\Pi_{U_\Box}$ be the decomposition subgroup corresponding to $J_\Box$.
The normalization of $J_\Box$ is isomorphic to $\mathbb{P}_{k_\Box}^1\setminus\{0, \infty\}$ if and only if the number of $k_\Box$-rational points of the normalization of $J_\Box$ is $q_{k_\Box}-1$.
Whether the latter condition holds or not is determined from $\Pi_{J_\Box}$.
\qed

\

\

\section{Fundamental groups of non-projective singular curves over generalized sub-$p$-adic fields}\label{GENP}

For $\Box\in\{\circ, \bullet, \emptyset\}$, we maintain the notation introduced at the beginning of \S 3.
Moreover, suppose that $k_\circ=k_\bullet$, and that $k_\circ=k_\bullet$ is a generalized sub-$p$-adic field (for some $p\in\mathfrak{Primes}$).

In the following, we suppose that $g_\Box+r_\Box-n_{\Box, \mathrm{a}}>0$ for $\Box\in\{\circ, \bullet, \emptyset\}$.
Suppose, moreover, that we are given an isomorphism $\phi:\Pi_{U_\circ}\stackrel{\sim}{\to}\Pi_{U_\bullet}$ of profinite groups which induces an isomorphism $\Delta_{U_\circ}\stackrel{\sim}{\to}\Delta_{U_\bullet}$ and the identity automorphism of $G:=G_{k_\circ}=G_{k_\bullet}$.
Then, by Remark \ref{p}, it holds that $\Sigma_\circ=\Sigma_\bullet$ (note that $k_\circ=k_\bullet$ is of characteristic zero).
Furthermore, suppose that $\Sigma_\circ=\Sigma_\bullet$ is an infinite set containing $p$.

\begin{Prop}\label{invp}
The following hold:
\begin{enumerate}
\item[(i)] $g_\circ=g_\bullet$.

\item[(ii)] $n^{>0}_{\circ, \mathrm{p}}=n^{>0}_{\bullet, \mathrm{p}}$.

\item[(iii)] $n^{>0}_{\circ, \mathrm{a}}=n^{>0}_{\bullet, \mathrm{a}}$.

\item[(iv)] $n^{0}_{\circ, >2}=n^{>0}_{\bullet, >2}$.

\item[(v)] $n^0_{\circ, 2}=n^0_{\bullet, 2}$.

\item[(vi)] $r_\circ-n^0_{\circ, 1}=r_\bullet-n^0_{\bullet, 1}$.

\item[(vii)] $\delta_\circ=\delta_\bullet$.
\end{enumerate}
\end{Prop}

\prf

Since we are given an isomorphism $\phi:\Pi_{U_\circ}\stackrel{\sim}{\to}\Pi_{U_\bullet}$ satisfying the above conditions, by considering the Hodge-Tate decomposition of the maximal pro-$p$ abelian quotient of $\Delta_{U_\Box}$ for $\Box\in\{\circ, \bullet\}$ (note that $p\in\Sigma_\circ=\Sigma_\bullet$), we obtain $g_\circ+\delta_\circ=g_\bullet+\delta_\bullet$ and $g_\circ+r_\circ-n_{\circ, \mathrm{a}}=g_\bullet+r_\bullet-n_{\bullet, \mathrm{a}}$ (cf. \cite{De}, \cite{Du}, and \cite[\S 4]{sur}).
First, we claim that the following assertion holds:

\begin{quote}
\noindent
\underline{Claim}

Suppose that $k$ is a generalized sub-$p$-adic field and $\Sigma$ is an infinite set containing $p$.
Then the kernel of the natural surjection $\Pi_U\twoheadrightarrow\Pi_X$ is recovered from $\Delta_U\subset\Pi_U$, and a family of invariants $\{g_\mathcal{H}+\delta_{\mathcal{H}}\}_{\mathcal{H}\subset\Pi_U}$, where $\mathcal{H}\subset\Pi_U$ runs through the set of open normal subgroups of $\Pi_U$.
\end{quote}
We shall prove the above claim.
Consider the isomorphism $\displaystyle\iota:\Delta_U\stackrel{\sim}{\to}\Delta=\left(\Star_{i=1}^n\Delta_i\right)\ast F_\delta$ (cf. (\ref{tisom})).
Let $\mathcal{H}$ be an open normal subgroup of $\Pi_U$, and set $H:=\iota(\Delta_U\cap\mathcal{H})$ and $N:=(\Delta_U:\Delta_U\cap\mathcal{H})$.
Then, by the Hurwitz formula, we obtain the following:
\begin{align*}
g_{\mathcal{H}}&=\sum_{i=1}^n|\Delta/H\Delta_i|(|H\Delta_i/H|(g_i-1)+1)+\sum_{x\in X_{\Delta_U\cap\mathcal{H}}}\dfrac{1}{2}(e_x-1) \\
&=\sum_{i=1}^nN(g_i-1)+\sum_{i=1}^n|\Delta/H\Delta_i|+\sum_{x\in X_{\Delta_U\cap\mathcal{H}}}\dfrac{1}{2}(e_x-1), \\
\delta_{\mathcal{H}}&=|\Delta/H\Delta_\delta|(|H\Delta_\delta/H|(\delta-1)+1)+1+nN-\sum_{i=1}^n|\Delta/H\Delta_i|-|\Delta/H\Delta_\delta|,
\end{align*}
where $e_x$ is the ramification index of the covering $X_{\Delta_U\cap\mathcal{H}}\to X_{\Delta_U}(=\overline{X})$ at $x\in X_{\Delta_U\cap\mathcal{H}}$.
Therefore,
\begin{align*}
g_{\mathcal{H}}+\delta_{\mathcal{H}}&=\sum_{i=1}^nN(g_i-1)+N(\delta-1+n)+1+\sum_{x\in X_{\Delta_U\cap\mathcal{H}}}\dfrac{1}{2}(e_x-1) \\
&=N(g+\delta-1)+1+\sum_{x\in X_{\Delta_U\cap\mathcal{H}}}\dfrac{1}{2}(e_x-1).
\end{align*}
For $m\in N(\Sigma)$, set:

\begin{align*}
\mathfrak{H}_m&:=\{\mathcal{H}\subset\Pi_U:\text{ an open normal subgroup}\,|\,(\Delta_U:\Delta_U\cap\mathcal{H})=m\}, \\
M_m&:=\min\{g_\mathcal{H}+\delta_\mathcal{H}\,|\,\mathcal{H}\in\mathfrak{H}_m\}, \\
\mathfrak{H}'_m&:=\{\mathcal{H}\in\mathfrak{H}\,|\,g_\mathcal{H}+\delta_{\mathcal{H}}=M_m\}, \\
\end{align*}
Note that these are recovered from $\Delta_U\subset\Pi_U$ and $\{g_{\mathcal{H}}+\delta_{\mathcal{H}}\}_{\mathcal{H}\subset\Pi_U}$.
Moreover, for $m\in N(\Sigma)$ and $\mathcal{H}\in\mathfrak{H}_m$, $\mathcal{H}\in\mathfrak{H}'_m$ if and only if the covering $X_{\Delta_U\cap\mathcal{H}}\to X_{\Delta_U}$ is unramified.
This shows that
\[\mathrm{Ker}(\Pi_U\twoheadrightarrow\Pi_X)=\mathrm{Ker}(\Delta_U\twoheadrightarrow\Delta_X)=\bigcap_{m\in N(\Sigma)}\bigcap_{\mathcal{H}\in\mathfrak{H}'_m}\mathcal{H},\]
as desired.

Therefore, the isomorphism $\phi:\Pi_{U_\circ}\stackrel{\sim}{\to}\Pi_{U_\bullet}$ induces an isomorphism $\Pi_{X_\circ}\stackrel{\sim}{\to}\Pi_{X_\bullet}$ (over $G$).
By considering the Hodge-Tate decomposition of the maximal pro-$p$ quotient of $\Delta_{X_\Box}$ for $\Box\in\{\circ, \bullet\}$, we obtain $g_\circ+\delta_\circ=g_\bullet+\delta_\bullet$ and $g_\circ=g_\bullet$, and hence also $\delta_\circ=\delta_\bullet$ and $r_\circ-n_{\circ, \mathrm{a}}=r_\bullet-n_{\bullet, \mathrm{a}}$.
This shows (i) and (vii), and (ii)-(vi) follow from similar arguments to the proof of Proposition \ref{inv}.
\qed

\begin{Cor}
The isomorphism $\phi|_{\Delta_{U_\circ}}:\Delta_{U_\circ}\stackrel{\sim}{\to}\Delta_{U_\bullet}$ (resp. $\iota_\bullet\circ\phi|_{\Delta_{U_\circ}}\circ\iota^{-1}_\circ:\Delta_\circ\stackrel{\sim}{\to}\Delta_\bullet$) induces an isomorphism $\psi:\Delta_{X_{\circ}}\stackrel{\sim}{\to}\Delta_{X_\bullet}$ (resp. $\iota_\delta:\Delta_{\delta_\circ}\stackrel{\sim}{\to}\Delta_{\delta_\bullet}$), i.e., the following diagrams are commutative:
\[\xymatrix@M=15pt{\Delta_{U_\circ} \ar[r]^{\sim}_{\phi|_{\Delta_{U_\circ}}} \ar@{->>}[d] & \Delta_{U_\bullet} \ar@{->>}[d] & & \Delta_\circ \ar[r]^{\sim}_{\iota_\bullet\circ\phi|_{\Delta_{U_\circ}}\circ\iota^{-1}_\circ} \ar@{->>}[d] & \Delta_\bullet \ar@{->>}[d] \\ \Delta_{X_\circ} \ar[r]^{\sim}_{\psi} & \Delta_{X_\bullet} & & \Delta_{\delta_\circ} \ar[r]^{\sim}_{\iota_\delta} & \Delta_{\delta_\bullet},}\]
where the vertical arrows are the natural surjections.
\end{Cor}

\prf

Immediate from Proposition \ref{invp}, its proof, and similar arguments to the proof of Corollary \ref{surj}.
\qed

\

In the following, for $\Box\in\{\circ, \bullet\}$, we assume that $(g_{\Box, i}, r_{\Box, i})\not\in\{(0, 0), (0, 1)\}$ for $1\leq i\leq n_\Box$ (i.e., $n^0_{\Box, 1}=n^0_{\Box, 0}=0$).

For $\Box\in\{\circ, \bullet\}$ and any open normal subgroup $\mathcal{H}_\Box\subset\Pi_{U_\Box}$, let $I_{U_{\mathcal{H}_\Box}}$ be the set of irreducible component of $U_{\mathcal{H}_\Box}$, and set:
\[\mathfrak{I}_\Box:=\lim_{\mathcal{H}_\Box\subset\Pi_{U_\Box}}I_{U_{\mathcal{H}_\Box}},\]
where $\mathcal{H}_\Box$ runs through the set of open normal subgroups of $\Pi_{U_\Box}$.
$\Pi_{U_\Box}$ acts naturally on $\mathfrak{I}_\Box$.

\begin{Prop}\label{irredp}
The isomorphism $\phi:\Pi_{U_\circ}\stackrel{\sim}{\to}\Pi_{U_\bullet}$ induces a bijection $\phi_\mathfrak{I}:\mathfrak{I}_\circ\to\mathfrak{I}_\bullet$ which is compatible with the natural actions of $\Pi_{U_\circ}$ and $\Pi_{U_\bullet}$ on $\mathfrak{I}_\circ$ and $\mathfrak{I}_\bullet$ (i.e., for any $\sigma\in\Pi_{U_\circ}$ and $J\in\mathfrak{I}_\circ$, it holds that $\phi_\mathfrak{I}(\sigma\cdot J)=\phi(\sigma)\cdot\phi_\mathfrak{I}(J)$).
\end{Prop}

\prf

Immediate from Proposition \ref{invp} and similar arguments to the proofs of Propositions \ref{irred}, \ref{et} and \ref{irredfin}.
\qed

\

In the following, for $\Box\in\{\circ, \bullet\}$ and $1\leq i\leq n_\Box$, we suppose that $\overline{U_{\Box, i}}^\nu$ is a hyperbolic curve (i.e., $2g_{\Box, i}+r_{\Box, i}-2>0$).

\begin{Thm}\label{mainsubp}
Set $k:=k_\circ=k_\bullet$.
Suppose that $k$ is a generalized sub-$p$-adic field, and that there exists an isomorphism $\phi:\Pi_{U_\circ}\stackrel{\sim}{\to}\Pi_{U_\bullet}$ which induces an isomorphism $\Delta_{U_\circ}\stackrel{\sim}{\to}\Delta_{U_\bullet}$ and the identity automorphism of the absolute Galois group $G_k$ of $k$.
Suppose, moreover, that $\Sigma_\circ=\Sigma_\bullet$ (cf. the beginning of the present section) is an infinite set containing $p$.
Then $U_\circ^\nu$ and $U_\bullet^\nu$ are isomorphic as $k$-schemes.
\end{Thm}

\prf

Immediate from \cite[Theorem 4.12]{sur}, Remark \ref{afdecomp} and Proposition \ref{irredp} (see also the proof of Theorem \ref{mainfin}).
\qed

\end{document}